\documentclass{article}

\usepackage{amsmath}
\usepackage{amsthm}
\usepackage{alltt}
\usepackage{hyperref}

\numberwithin{equation}{section}

\newcommand{\res}[1]
{\textrm{\rm Res}_{#1}}

\title{Some Weighted Generalized Fibonacci Number Summation Identities, Part 2}

\author{M.J. Kronenburg}
\date{}

\begin{document}

\maketitle

\begin{abstract}
In part 1 of this paper some  linear weighted generalized Fibonacci number summation identities were
derived using the fact that the Fibonacci number is the residue of a rational function.
In this part, using the same method, some quadratic and cubic weighted generalized Fibonacci number summation
identities are derived, including some infinite series and generating functions.
In addition some quadratic and cubic weighted generalized Fibonacci number summation identities
with binomial coefficients, and some cubic and quartic Fibonacci and Lucas number summation identities are derived.
\end{abstract}


\section{Quadratic weighted generalized Fibonacci\\ number summation identities}

Let $F_n$ be the Fibonacci number, $L_n$ be the Lucas number, and let $G_n$ be a generalized Fibonacci number
for which $G_{n+2}=G_{n+1}+G_n$ with any seed $G_0$ and $G_1$,
and let $H_n$ also be a generalized Fibonacci number for which $H_{n+2}=H_{n+1}+H_n$
with any seed $H_0$ and $H_1$ \cite{BQ03,CW00,V89}.
The following are quadratic weighted generalized Fibonacci number summation identities.
\begin{equation}\label{delta1}
 \Delta_1 = 1-L_{b+d} x + (-1)^{b+d} x^2
\end{equation}
\begin{equation}\label{delta2}
 \Delta_2 = 1-(-1)^d L_{b-d} x + (-1)^{b-d} x^2
\end{equation}
\begin{equation}\label{pres}
 P_2(v,w) = \sum_{k=0}^v \binom{v}{k} (-1)^{(b+d+1)k} x^k G_{a+bw-bk}H_{c+dw-dk}
\end{equation}
\begin{equation}\label{qres}
 Q_2(v,w) = \frac{1}{5} [ ( \frac{2}{\Delta_1^v} + \frac{3}{\Delta_2^v} ) P_2(v,w) + ( \frac{1}{\Delta_1^v} - \frac{1}{\Delta_2^v} ) ( P_2^-(v,w) + P_2^+(v,w) ) ]
\end{equation}
where $P_2^-(v,w)$ resp. $P_2^+(v,w)$ is identical to $P_2(v,w)$ but with
$a$ replaced by $a-1$ and $c$ by $c-1$ resp. $a$ by $a+1$ and $c$ by $c+1$.
\begin{equation}\label{res1}
 \sum_{k=0}^n x^k G_{a+bk}H_{c+dk} = - x^{n+1} Q_2(1,n+1) + Q_2(1,0)
\end{equation}
\begin{equation}\label{res2}
 \sum_{k=0}^n x^k k G_{a+bk}H_{c+dk} = -(n+1)x^{n+1}Q_2(1,n+1) - [ x^{n+2}Q_2(2,n+2) - xQ_2(2,1) ]
\end{equation}
\begin{equation}\label{res3}
\begin{split}
 \sum_{k=0}^n x^k k^2 G_{a+bk}H_{c+dk} = & -(n+1)^2 x^{n+1} Q_2(1,n+1) \\
 & - [ (2n+3) x^{n+2} Q_2(2,n+2) - x Q_2(2,1) ] \\
 & - 2 [ x^{n+3} Q_2(3,n+3) - x^2 Q_2(3,2) ] \\
\end{split}
\end{equation}
\begin{equation}\label{res4}
\begin{split}
 \sum_{k=0}^n x^k k^3 G_{a+bk}H_{c+dk} = & -(n+1)^3 x^{n+1} Q_2(1,n+1) \\
 & - [ (3n^2+9n+7) x^{n+2} Q_2(2,n+2) - x Q_2(2,1) ] \\
 & - 6 [ (n+2) x^{n+3} Q_2(3,n+3) - x^2 Q_2(3,2) ] \\
 & - 6 [ x^{n+4} Q_2(4,n+4) - x^3 Q_2(4,3) ] \\ 
\end{split}
\end{equation}
\begin{equation}\label{res5}
\begin{split}
 \sum_{k=0}^n & x^k k^4 G_{a+bk}H_{c+dk} \\
= & -(n+1)^4 x^{n+1} Q_2(1,n+1) \\
 & - [ (4n^3+18n^2+28n+15) x^{n+2} Q_2(2,n+2) - x Q_2(2,1) ] \\
 & - 2 [ (6n^2+24n+25) x^{n+3} Q_2(3,n+3) - 7 x^2 Q_2(3,2) ] \\
 & - 12 [ (2n+5) x^{n+4} Q_2(4,n+4) - 3 x^3 Q_2(4,3) ] \\
 & - 24 [ x^{n+5} Q_2(5,n+5) - x^4 Q_2(5,4) ] \\ 
\end{split}
\end{equation}
\begin{equation}\label{res6}
\begin{split}
 \sum_{k=0}^n & x^k k^5 G_{a+bk}H_{c+dk} \\
= & -(n+1)^5 x^{n+1} Q_2(1,n+1) \\
 & - [ (5n^4+30n^3+70n^2+75n+31) x^{n+2} Q_2(2,n+2) - x Q_2(2,1) ] \\
 & - 10 [ (2n^3+12n^2+25n+18) x^{n+3} Q_2(3,n+3) - 3 x^2 Q_2(3,2) ] \\
 & - 30 [ (2n^2+10n+13) x^{n+4} Q_2(4,n+4) - 5 x^3 Q_2(4,3) ] \\
 & - 120 [ (n+3) x^{n+5} Q_2(5,n+5) - 2 x^4 Q_2(5,4) ] \\
 & - 120 [ x^{n+6} Q_2(6,n+6) - x^5 Q_2(6,5) ] \\
\end{split}
\end{equation}
When $d=0$ then the formulas of part 1 of this paper \cite{K19} reappear
because then $\Delta_1=\Delta_2=\Delta$, $P_2(v,w)=H_c P(v,w)$ and $Q_2(v,w)=H_c\Delta^{-v}P(v,w)$,
and likewise when $b=0$ because $(-1)^d L_{-d} = L_d$.\\

\section{The Special Cases $d=b$ and $d=-b$}

In the special cases $d=b$ and $d=-b$ the formula for $Q_2(v,w)$ takes a special form.\\
The special case $d=b$:
\begin{equation}
 \Delta_1 = 1 - L_{2b} x + x^2
\end{equation}
\begin{equation}
 \Delta_2 = ( 1 - (-1)^b x )^2
\end{equation}
\begin{equation}
 P_2(v,w) = \sum_{k=0}^v \binom{v}{k} (-1)^k x^k G_{a+bw-bk}H_{c+bw-bk}
\end{equation}
\begin{equation}\label{quadspec1}
\begin{split}
 & Q_2(v,w) = \frac{1}{\Delta_1^v} P_2(v,w) \\
 & + \frac{1}{5} ( \frac{1}{\Delta_2^v} - \frac{1}{\Delta_1^v} ) (-1)^{bw} ( 1 - (-1)^b x )^v ( 3G_aH_c - G_{a-1}H_{c-1} - G_{a+1}H_{c+1} ) \\
\end{split}
\end{equation}
The special case $d=-b$:
\begin{equation}
 \Delta_1 = ( 1 - x )^2
\end{equation}
\begin{equation}
 \Delta_2 = 1 - (-1)^b L_{2b} x + x^2
\end{equation}
\begin{equation}
 P_2(v,w) = \sum_{k=0}^v \binom{v}{k} (-1)^k x^k G_{a+bw-bk}H_{c-bw+bk}
\end{equation}
\begin{equation}\label{quadspec2}
  Q_2(v,w) = \frac{1}{\Delta_2^v} P_2(v,w) 
  + \frac{1}{5} ( \frac{1}{\Delta_1^v} - \frac{1}{\Delta_2^v} ) ( 1 - x )^v ( 2G_aH_c + G_{a-1}H_{c-1} + G_{a+1}H_{c+1} )
\end{equation}
When $d=b$ and $x=(-1)^b$ or $d=-b$ and $x=1$ then $\Delta_2$ or $\Delta_1$ is zero.
In that case the corresponding limit $x\to (-1)^b$ or $x\to 1$
must be taken, applying the following formula which follows from the Taylor expansion of $f(x)$,
where $f^{(n)}(x)$ is the $n$-th derivative of $f(x)$:
\begin{equation}\label{taylor}
 \frac{1}{n!} f^{(n)}(x_0) = \lim_{x\to x_0} \frac{1}{(x-x_0)^n} [ f(x) - \sum_{k=0}^{n-1} \frac{1}{k!} (x-x_0)^k f^{(k)}(x_0) ]
\end{equation}
For $x\neq -1$:
\begin{equation}
\begin{split}
 & \sum_{k=0}^n x^k G_{p+k}H_{q+k} \\
  = & \frac{1}{1-3x+x^2} [ x^{n+1} ( x G_{p+n}H_{q+n} - G_{p+n+1}H_{q+n+1} ) + G_pH_q - x G_{p-1}H_{q-1} ] \\
 & + \frac{x [ (-1)^{n+1} x^{n+1} - 1 ]}{(1+x)(1-3x+x^2)} ( 3G_pH_q - G_{p-1}H_{q-1} - G_{p+1}H_{q+1} ) \\
\end{split}
\end{equation}
For $x=-1$ taking $f(x)=(-1)^{n+1}x^{n+1}$, $x_0=-1$ and $f^{(1)}(-1)=-(n+1)$ the result is (\ref{reslist1}).\\
For $x\neq 1$:
\begin{equation}
\begin{split}
 & \sum_{k=0}^n x^k G_{p+k}H_{q-k} \\
  = & \frac{1}{1+3x+x^2} [ x^{n+1} ( x G_{p+n}H_{q-n} - G_{p+n+1}H_{q-n-1} ) + G_pH_q - x G_{p-1}H_{q+1} ] \\
 & - \frac{x ( x^{n+1} - 1 )}{(1-x)(1+3x+x^2)} ( 2G_pH_q + G_{p-1}H_{q-1} + G_{p+1}H_{q+1} ) \\
\end{split}
\end{equation}
For $x=1$ taking $f(x)=x^{n+1}$, $x_0=1$ and $f^{(1)}(1)=n+1$ the result is (\ref{resultconv1}).\\

\section{Cubic weighted generalized Fibonacci number\\ summation identities}

Let $G_n$ be a generalized Fibonacci number
for which $G_{n+2}=G_{n+1}+G_n$ with any seed $G_0$ and $G_1$,
and let $H_n$ also be a generalized Fibonacci number for which $H_{n+2}=H_{n+1}+H_n$
with any seed $H_0$ and $H_1$,
and let $K_n$ also be a generalized Fibonacci number for which $K_{n+2}=K_{n+1}+K_n$
with any seed $K_0$ and $K_1$  \cite{BQ03,CW00,V89}.
The following are cubic weighted generalized Fibonacci number summation identities.
\begin{equation}\label{cubicres1}
\begin{split}
 \sum_{k=0}^n & x^k G_{p+k}H_{q+k}K_{r+k} = \frac{1}{(1+x-x^2)(1-4x-x^2)} \\
 \cdot & \{ x^{n+1} [ x^2(x+3) G_{p+n}H_{q+n}K_{r+n} + (3x-1) G_{p+n+1}H_{q+n+1}K_{r+n+1} \\
 & \qquad + x^2 G_{p+n-1}H_{q+n-1}K_{r+n-1} - x G_{p+n+2}H_{q+n+2}K_{r+n+2} ] \\
 & - x^2(x+3) G_{p-1}H_{q-1}K_{r-1} - (3x-1) G_pH_qK_r \\
 & - x^2 G_{p-2}H_{q-2}K_{r-2} + x G_{p+1}H_{q+1}K_{r+1} \} \\
\end{split}
\end{equation}
\begin{equation}\label{cubicres2}
\begin{split}
 \sum_{k=0}^n & x^k G_{p+k}H_{q+k}K_{r-k} = \frac{1}{(1-x-x^2)(1+4x-x^2)} \\
 \cdot & \{ x^{n+1} [ x^2(x-3) G_{p+n}H_{q+n}K_{r-n} - (3x+1) G_{p+n+1}H_{q+n+1}K_{r-n-1} \\
 & \qquad + x^2 G_{p+n-1}H_{q+n-1}K_{r-n+1} - x G_{p+n+2}H_{q+n+2}K_{r-n-2} ] \\
 & - x^2(x-3) G_{p-1}H_{q-1}K_{r+1} + (3x+1) G_pH_qK_r \\
 & - x^2 G_{p-2}H_{q-2}K_{r+2} + x G_{p+1}H_{q+1}K_{r-1} \} \\
\end{split}
\end{equation}
\begin{equation}\label{cubicres3}
\begin{split}
 \sum_{k=0}^n & x^k G_{p+k}H_{q-k}K_{r-k} = \frac{1}{(1+x-x^2)(1-4x-x^2)} \\
 \cdot & \{ x^{n+1} [ x^2(x+3) G_{p+n}H_{q-n}K_{r-n} + (3x-1) G_{p+n+1}H_{q-n-1}K_{r-n-1} \\
 & \qquad + x^2 G_{p+n-1}H_{q-n+1}K_{r-n+1} - x G_{p+n+2}H_{q-n-2}K_{r-n-2} ] \\
 & - x^2(x+3) G_{p-1}H_{q+1}K_{r+1} - (3x-1) G_pH_qK_r \\
 & - x^2 G_{p-2}H_{q+2}K_{r+2} + x G_{p+1}H_{q-1}K_{r-1} \} \\
\end{split}
\end{equation}
\begin{equation}\label{cubicres4}
\begin{split}
 & \sum_{k=0}^n x^k G_{p+2k}H_{q+2k}K_{r+2k} = \frac{1}{(1-3x+x^2)(1-18x+x^2)} \\
 \cdot & \{ x^{n+1} [ x^2(x-21) G_{p+2n}H_{q+2n}K_{r+2n} + (21x-1) G_{p+2n+2}H_{q+2n+2}K_{r+2n+2} \\
 & \qquad + x^2 G_{p+2n-2}H_{q+2n-2}K_{r+2n-2} - x G_{p+2n+4}H_{q+2n+4}K_{r+2n+4} ] \\
 & - x^2(x-21) G_{p-2}H_{q-2}K_{r-2} - (21x-1) G_pH_qK_r \\
 & - x^2 G_{p-4}H_{q-4}K_{r-4} + x G_{p+2}H_{q+2}K_{r+2} \} \\
\end{split}
\end{equation}

\section{Generalized Fibonacci Identities}

The following is a list of generalized Fibonacci identities 
some of which are used in the derivation.

\begin{equation}\label{gfl}
 G_n = \frac{1}{2} [ ( 2G_1 - G_0 ) F_n + G_0 L_n ]
\end{equation}
\begin{equation}\label{gmin}
 G_{-n} = (-1)^{n+1} \frac{1}{2} [ ( 2G_1 - G_0 ) F_n - G_0 L_n ]
\end{equation}
\begin{equation}\label{gfromf}
 G_n = ( G_1 - G_0 ) F_n + G_0 F_{n+1}
\end{equation}
\begin{equation}\label{gfroml}
 G_n = \frac{1}{5} [ (3G_0-G_1) L_n + (2G_1-G_0) L_{n+1} ]
\end{equation}
\begin{equation}\label{ffromg}
 F_n = \frac{ -G_1G_n + G_0G_{n+1} }{ G_0G_2 - G_1^2 }
\end{equation}
\begin{equation}\label{lfromg}
 L_n = \frac{ (2G_0+G_1) G_n + (G_0-2G_1) G_{n+1} }{ G_0G_2 - G_1^2 }
\end{equation}
\begin{equation}\label{gminus}
 G_{-n} = (-1)^n \frac{(G_0^2+G_1^2)G_n + (G_{-1}^2-G_1^2)G_{n+1}}{G_0G_2 - G_1^2}
\end{equation}
\begin{equation}\label{hfromg}
 H_n = \frac{(H_0G_2 - H_1G_1) G_n + (G_0H_1 - G_1H_0) G_{n+1}}{G_0G_2 - G_1^2}
\end{equation}
\begin{equation}\label{ghprod}
\begin{split}
 G_nH_m = \frac{1}{5} \{ & (3G_0-G_1)H_{m+n} + (2G_1-G_0)H_{m+n+1} \\
 & + (-1)^n [ (2G_0+G_1)H_{m-n} - (2G_1-G_0)H_{m-n+1} ] \} \\
\end{split}
\end{equation}
\begin{equation}
 G_n^2 = \frac{1}{5} [ (3G_0-G_1) G_{2n} + (2G_1-G_0) G_{2n+1} + 2 (-1)^n (G_0G_2-G_1^2) ]
\end{equation}
\begin{equation}\label{comm}
 G_nH_m - G_mH_n = (-1)^n ( G_0H_1 - G_1H_0 ) F_{m-n}
\end{equation}
\begin{equation}\label{general}
 G_{n+r+p}H_{m+r+q} - G_{n+r}H_{m+r+p+q} = (-1)^{n+r} ( G_pH_{m-n+q} - G_0H_{m-n+p+q} )
\end{equation}
\begin{equation}\label{plusplus}
 G_{n-p}H_{m+p} + G_{n-p+1}H_{m+p+1} = G_0H_{m+n} + G_1H_{m+n+1}
\end{equation}
\begin{equation}\label{cassini}
 G_{n+p}H_{m+p} - G_{n+p-1}H_{m+p+1} = (-1)^{n+p} ( G_0H_{m-n} - G_{-1}H_{m-n+1} )
\end{equation}
\begin{equation}\label{prod1}
\begin{split}
 G_nH_m + G_{n+1}H_{m+1} = \frac{1}{5} & [ (3G_0-G_1)(H_{m+n}+H_{m+n+2}) \\
  & + (2G_1-G_0)(H_{m+n+1}+H_{m+n+3}) ] \\
\end{split}
\end{equation}
\begin{equation}\label{prod2}
\begin{split}
 G_nH_m - G_{n-1}H_{m+1} = (-1)^n \frac{1}{5} & [ (2G_0+G_1)(H_{m-n}+H_{m-n+2}) \\
  & - (2G_1-G_0)(H_{m-n+1}+H_{m-n+3}) ] \\
\end{split}
\end{equation}
\begin{equation}\label{prod3}
\begin{split}
 2G_nH_m + G_{n-1}H_{m-1} + & G_{n+1} H_{m+1} = 3G_nH_m + G_{n-1}H_{m+1} + G_{n+1}H_{m-1} \\
 = & (3G_0-G_1) H_{m+n} + (2G_1-G_0) H_{m+n+1} \\
\end{split}
\end{equation}
\begin{equation}\label{prod4}
\begin{split}
 3G_nH_m - G_{n-1}H_{m-1} - & G_{n+1}H_{m+1} = 2G_nH_m - G_{n-1}H_{m+1} - G_{n+1}H_{m-1} \\
 = & (-1)^n [ (2G_0+G_1) H_{m-n} - (2G_1-G_0) H_{m-n+1} ] \\
\end{split}
\end{equation}
The identity (\ref{gfl}) is easily proved by checking it for $n=0$ and $n=1$ and $G_{n+2}=G_{n+1}+G_n$.
The identity (\ref{gmin}) follows from (\ref{gfl}) and $F_{-n}=(-1)^{n+1}F_n$ and $L_{-n}=(-1)^nL_n$.
The identities (\ref{gfromf}) and (\ref{gfroml}) follow from (\ref{gfl}) and $L_n=2F_{n+1}-F_n$
and $5F_n=2L_{n+1}-L_n$.
These two identities form a 2x2 matrix equation, and inverting this matrix yields (\ref{ffromg}) and (\ref{lfromg}).
Substituting these two identities in (\ref{gmin}) yields (\ref{gminus}).
The identity (\ref{hfromg}) follows from (\ref{gfl}) with $G=H$, (\ref{ffromg}) and (\ref{lfromg}).
The identity (\ref{ghprod}) with $H=G$ is derived below in the text of this paper,
and is generalized using (\ref{hfromg}).
The identity (\ref{comm}) is derived for $H=F$ and $H=L$ by (\ref{ghprod}) with $F_{-n}=(-1)^{n+1}F_n$ and $L_{-n}=(-1)^nL_n$,
and then (\ref{gfl}) with $G=H$ is used.
The identities (\ref{general}) to (\ref{prod2}) are directly derived from (\ref{ghprod}).
The identities (\ref{prod3}) and (\ref{prod4}) are also derived from (\ref{ghprod}) 
using (\ref{prod1}) and (\ref{prod2}) and using $H_{n-2}+2H_n+H_{n+2}=5H_n$.\\
In literature \cite{V89}, (21) is the special case $H=F$ of (\ref{comm}),
(18) is the special case $m=n$ and $r=0$ of (\ref{general}),
and (28) is the special case $H=G$, $m=n$ and $p=0$ of (\ref{cassini}).

\section{Derivation of the Quadratic Summation\\ Identities}

The following sum is to be evaluated:
\begin{equation}
\begin{split}
 \sum_{k=0}^n & x^k F_{a+bk}F_{c+dk} = \sum_{k=0}^n x^k \res{y}\frac{(1+y)^{a+bk}}{(y+\phi)(y-1/\phi)}\res{z}\frac{(1+z)^{c+dk}}{(z+\phi)(z-1/\phi)} \\
 & = \res{y}\res{z} \frac{(1+y)^a(1+z)^c}{(y+\phi)(y-1/\phi)(z+\phi)(z-1/\phi)} \sum_{k=0}^n [x(1+y)^b(1+z)^d]^k \\
 & = \res{y}\res{z} \frac{(1+y)^a(1+z)^c}{(y+\phi)(y-1/\phi)(z+\phi)(z-1/\phi)} \frac{[x(1+y)^b(1+z)^d]^{n+1}-1}{x(1+y)^b(1+z)^d-1} \\
\end{split}
\end{equation}
Taking the residues $y=z=1/\phi$ and $y=z=-\phi$, and using theorem 5.1 in \cite{K19} with $v=-1$, $w=x$ and $n=b+d$,
the result is:
\begin{equation}
\begin{split}
 Z_1 = \frac{1}{5\Delta_1} & \{ x^{n+1} [ (-1)^{b+d} L_{a+c+(b+d)n}x - L_{a+c+(b+d)(n+1)} ] \\
 & - [ (-1)^{b+d} L_{a+c-(b+d)} - L_{a+c} ] \} \\
\end{split}
\end{equation}
and taking the residues $y=1/\phi,z=-\phi$ and $y=-\phi,z=1/\phi$, 
and using $\phi^a(1-\phi)^c=(-1)^c\phi^{a-c}$ and $(1-\phi)^a\phi^c=(-1)^c(1-\phi)^{a-c}$
and $\phi^b(1-\phi)^d=(-1)^d\phi^{b-d}$ and $(1-\phi)^b\phi^d=(-1)^d(1-\phi)^{b-d}$,
and using theorem 5.1 in \cite{K19} with $v=-1$, $w=(-1)^d x$ and $n=b-d$,
the result is:
\begin{equation}
\begin{split}
 Z_2 = -\frac{1}{5\Delta_2} (-1)^c & \{ x^{n+1} (-1)^{d(n+1)} [ (-1)^b L_{a-c+(b-d)n}x - L_{a-c+(b-d)(n+1)} ] \\
 & - [ (-1)^b L_{a-c-(b-d)} - L_{a-c} ] \} \\
\end{split}
\end{equation}
and the total result is:
\begin{equation}
 \sum_{k=0}^n x^k F_{a+bk}F_{c+dk} = Z_1 + Z_2
\end{equation}
This formula is now generalized for $G$ instead of $F$.
Let $\alpha=\frac{1}{2}(G_{-1}+G_1)$ and $\beta=\frac{1}{2}G_0$, then (\ref{gfl}) from \cite{K19} means:
\begin{equation}
 G_n = \alpha F_n + \beta L_n
\end{equation}
so that:
\begin{equation}
 G_m G_n = \alpha^2 F_m F_n + \beta^2 L_m L_n + \alpha\beta ( F_m L_n + L_m F_n )
\end{equation}
Using $L_n=F_{n-1}+F_{n+1}$ we have:
\begin{equation}
 L_m L_n = F_{m-1}F_{n-1} + F_{m+1}F_{n+1} + F_{m-1}F_{n+1} + F_{m+1}F_{n-1}
\end{equation}
Adding the corresponding $Z_1$ and $Z_2$ above, and using:
\begin{equation}
 L_{m+n-2} + L_{m+n+2} + 2 L_{m+n} = 5 L_{m+n}
\end{equation}
we obtain:
\begin{equation}
 \sum_{k=0}^n x^k L_{a+bk}L_{c+dk} = 5 ( Z_1 - Z_2 )
\end{equation}
Likewise we have:
\begin{equation}
 F_m L_n = F_{m-1} F_n + F_{m+1} F_n
\end{equation}
and using:
\begin{equation}
 L_{m+n-1} + L_{m+n+1} = 5 F_{m+n}
\end{equation}
and let $Z_1^*$ resp. $Z_2^*$ be $Z_1$ resp. $Z_2$ but with $L$ replaced by $F$,
then:
\begin{equation}
 \sum_{k=0}^n x^k F_{a+bk}L_{c+dk} = 5 ( Z_1^* - Z_2^* )
\end{equation}
\begin{equation}
 \sum_{k=0}^n x^k L_{a+bk}F_{c+dk} = 5 ( Z_1^* + Z_2^* )
\end{equation}
Combining these formulas it is clear that:
\begin{equation}\label{sumres}
 \sum_{k=0}^n x^k G_{a+bk}G_{c+dk} = \alpha^2 (Z_1 + Z_2) + 5\beta^2 (Z_1 - Z_2 ) + 10 \alpha\beta Z_1^*
\end{equation}
The right side can be expressed using $G$ instead of $L$ and $F$ using (\ref{ffromg}) and (\ref{lfromg}).
Since $\alpha^2-5\beta^2=G_{-1}G_1-G_0^2=G_1^2-G_0G_2$, we have:
\begin{equation}
 ( \alpha^2 - 5\beta^2 ) L_n = -(2G_0+G_1) G_n + (2G_1-G_0) G_{n+1}
\end{equation}
The other term evaluates to:
\begin{equation}
 ( \alpha^2 + 5\beta^2 ) L_n + 10\alpha\beta F_n = 10 \beta ( \beta L_n + \alpha F_n ) + ( \alpha^2 - 5\beta^2 ) L_n
\end{equation}
which is evaluated with:
\begin{equation}
 \beta L_n + \alpha F_n = ( G_{-1}G_1 - G_0^2 ) G_n
\end{equation}
so that:
\begin{equation}
 ( \alpha^2 + 5\beta^2 ) L_n + 10\alpha\beta F_n = ( 3G_0-G_1 ) G_n + ( 2G_1-G_0 ) G_{n+1} 
\end{equation}
Substituting these formulas in (\ref{sumres}) and taking $x=1$, $b=d=1$, $a=p$ and $c=q$
we obtain:
\begin{equation}\label{firstsum}
\begin{split}
 & \sum_{k=0}^n G_{p+k}G_{q+k} \\
= & \frac{1}{5} [ ( G_{p+q+2n+1} - G_{p+q-1} ) (3G_0-G_1) + ( G_{p+q+2n+2} - G_{p+q} ) (2G_1-G_0) ] \\
 & + \frac{1}{10} (-1)^q [(-1)^n+1] [ G_{p-q} (2G_0+G_1) - G_{p-q+1} (2G_1-G_0) ] \} \\
\end{split}
\end{equation}
The $n=0$ case of this formula and taking $p=m$ and $q=n$ yields the important identity:
\begin{equation}
\begin{split}
 G_nG_m = \frac{1}{5} \{ & (3G_0-G_1)G_{m+n} + (2G_1-G_0)G_{m+n+1} \\
 & + (-1)^n [ (2G_0+G_1)G_{m-n} - (2G_1-G_0)G_{m-n+1} ] \} \\
\end{split}
\end{equation}
This identity can be generalized to (\ref{ghprod}) using (\ref{hfromg}),
and interchanging $G$ and $H$:
\begin{equation}\label{prod}
\begin{split}
 G_mH_n = \frac{1}{5} \{ & (3H_0-H_1)G_{m+n} + (2H_1-H_0)G_{m+n+1} \\
 & + (-1)^n [ (2H_0+H_1)G_{m-n} - (2H_1-H_0)G_{m-n+1} ] \} \\
\end{split}
\end{equation}
Substituting this identity in the summand with $m=a+bk$ and $n=c+dk$,
the formulas of \cite{K19} can be applied
by replacing $a$ by $a+c$ or $a+c+1$, $b$ by $b+d$, or replacing $a$ by $a-c$ or $a-c+1$,
$b$ by $b-d$ and $x$ by $(-1)^d x$, and using the polynomials as derived in (5.26) in \cite{K19}:
\begin{equation}\label{r1}
 R_1(v,w) = \sum_{k=0}^v \binom{v}{k} (-1)^{(b+d+1)k} x^k G_{a+c+(b+d)(w-k)}
\end{equation}
\begin{equation}\label{r2}
 R_2(v,w) = \sum_{k=0}^v \binom{v}{k} (-1)^{(b+1)k} x^k G_{a-c+(b-d)(w-k)}
\end{equation}
\begin{equation}
 S_1(v,w) = (3H_0-H_1) R_1(v,w) + (2H_1-H_0) R_1^+(v,w)
\end{equation}
\begin{equation}
 S_2(v,w) = (2H_0+H_1) R_2(v,w) - (2H_1-H_0) R_2^+(v,w)
\end{equation}
where $R_1^+(v,w)$ resp. $R_2^+(v,w)$ is identical to $R_1(v,w)$ resp. $R_2(v,w)$
but with $a$ replaced by $a+1$.
The terms of the right sides of the quadratic identities have the form:
\begin{equation}\label{quadexpr}
 \Delta_1^{-v} x^w S_1(v,w) + \Delta_2^{-v} (-1)^{c+dw} x^w S_2(v,w)
\end{equation}
Using this formula and substituting the previous 4 formulas
and applying in its summand formulas (\ref{prod3}) and (\ref{prod4})
with $G$ and $H$ interchanged:
\begin{equation}\label{prodsec1}
 (3H_0-H_1) G_{m+n} + (2H_1-H_0) G_{m+n+1} = 2G_mH_n + G_{m-1}H_{n-1} + G_{m+1}H_{n+1}
\end{equation}
\begin{equation}\label{prodsec2}
 (-1)^n [ (2H_0+H_1) G_{m-n} - (2H_1-H_0) G_{m-n+1} ] = 3G_mH_n - G_{m-1}H_{n-1} - G_{m+1}H_{n+1} 
\end{equation}
with $m=a+bw-bk$ and $n=c+dw-dk$,
the formulas (\ref{pres}) to (\ref{res6}) have been derived.\\
For the special case $d=b$, taking the summands of (\ref{qres}) and using (\ref{prodsec2}):
\begin{equation}
\begin{split}
 & 3G_{a+bw-bk}H_{c+bw-bk} - G_{a+bw-bk-1}H_{c+bw-bk-1} - G_{a+bw-bk+1}H_{c+bw-bk+1} \\
 & = (-1)^{c+bw+bk} [ ( 2H_0+H_1 ) G_{a-c} - ( 2H_1-H_0 ) G_{a-c+1} ) \\
 & = (-1)^{bw+bk} ( 3G_aH_c - G_{a-1}H_{c-1} - G_{a+1}H_{c+1} ) \\
\end{split}
\end{equation}
\begin{equation}
\begin{split}
 & 2G_{a+bw-bk}H_{c+bw-bk} + G_{a+bw-bk-1}H_{c+bw-bk-1} + G_{a+bw-bk+1}H_{c+bw-bk+1} \\
 & = 5G_{a+bw-bk}H_{c+bw-bk} - ( 3G_{a+bw-bk}H_{c+bw-bk} \\
 & \qquad- G_{a+bw-bk-1}H_{c+bw-bk-1} - G_{a+bw-bk+1}H_{c+bw-bk+1} ) \\
 & = 5G_{a+bw-bk}H_{c+bw-bk} - (-1)^{bw+bk} ( 3G_aH_c - G_{a-1}H_{c-1} - G_{a+1}H_{c+1} ) \\
\end{split}
\end{equation}
\begin{equation}
 \sum_{k=0}^n \binom{v}{k} (-1)^{(b+1)k} x^k = (1-(-1)^bx)^v
\end{equation}
which yields (\ref{quadspec1}).
For the special case $d=-b$ a similar derivation using (\ref{prodsec1}) yields (\ref{quadspec2}).

\section{Derivation of the Cubic Summation Identities}

The cubic generalized Fibonacci summation identities can be derived
as in the previous section, by using again (\ref{prod})
and using a linear and a quadratic summation identity already derived.
The summand of:
\begin{equation}
 \sum_{k=0}^n x^k G_{a+k}H_{b+k}K_{c+k}
\end{equation}
is expressed with (\ref{prod}) into:
\begin{equation}
\begin{split}
 G_{a+k}H_{b+k}K_{c+k} & = \frac{1}{5} \{ [ (3H_0-H_1) G_{a+b+2k} + (2H_1-H_0) G_{a+b+2k+1} ] K_{c+k} \\
 & + (-1)^b [ (2H_0+H_1) G_{a-b} - (2H_1-H_0) G_{a-b+1} ] (-1)^k K_{c+k} \\
\end{split}
\end{equation}
The cubic sum can now be evaluated with a quadratic and a linear sum,
which after simplification become:
\begin{equation}
\begin{split}
 \sum_{k=0}^n & x^k G_{a+b+2k}K_{c+k} = \frac{1}{(1-4x-x^2)(1+x-x^2)} \\
 \cdot & \{ x^{n+1} [ x^3 G_{a+b+2n}K_{c+n} + x^2 ( G_{a+b+2n+2}K_{c+n} + G_{a+b+2n-2}K_{c+n+1} ) \\
 & \qquad + x ( G_{a+b+2n}K_{c+n+1} - G_{a+b+2n+4}K_{c+n} ) - G_{a+b+2n+2}K_{c+n+1} ] \\
 & - [ x^3 G_{a+b-2}K_{c-1} + x^2 ( G_{a+b}K_{c-1} + G_{a+b-4}K_c ) \\
 & \qquad + x ( G_{a+b-2}K_c - G_{a+b+2}K_{c-1} ) - G_{a+b}K_c ] \} \\
\end{split}
\end{equation}
\begin{equation}
 \sum_{k=0}^n x^k (-1)^k K_{c+k} = \frac{1}{1+x-x^2} [ (-1)^n x^{n+1} ( K_{c+n+1} - x K_{c+n} ) + K_c - x K_{c-1} ]
\end{equation}
\begin{equation}
\begin{split}
 & \sum_{k=0}^n x^k (-1)^k K_{c+k} = \frac{1}{(1-4x-x^2)(1+x-x^2)} \\
 \cdot & \{ (-1)^n x^{n+1} [ x^3 K_{c+n} + x^2 ( K_{c+n+2} - K_{c+n-3} )
   - x ( K_{c+n-1} + K_{c+n+4} ) + K_{c+n+1} ] \\
 & + x^3 K_{c-1} + x^2 ( K_{c+1} - K_{c-4} ) - x ( K_{c-2} + K_{c+3} ) + K_c \} \\
\end{split}
\end{equation}
For each power of $x$,
the terms with $G_{a+b+2n+2p}$ can now be expressed as a sum of three products 
of the form $G_{a+n+p}H_{b+n+p}$ using (\ref{prodsec1}),
and likewise for $G_{a+b+2p}$ as products $G_{a+p}H_{b+p}$.
The terms with $G_{a-b}$ can be expressed as a sum of three products
of the form $G_{a+n+p}H_{b+n+p}$ or $G_{a+p}H_{b+p}$ using (\ref{prodsec2}).
For each power of $x$, the resulting sum of triple products can be simplified to
the coefficient of that power of $x$ in (\ref{cubicres1}).

\section{Derivation of some Cubic and Quartic Fibonacci and Lucas Number Summation Identities}

Summation identities for $G=H=K=F$ or $G=H=K=L$ can be obtained by using the following identities derived
from (2.8), (2.9) and (2.10) of \cite{K18} or (15b), (17a) and (17b) of \cite{V89}
in the summand and applying summation identities with $G=F$ or $G=L$ from
part 1 of this paper \cite{K19}:
\begin{equation}\label{summand1}
 F_{p+k}F_{q+k} = \frac{1}{5} ( L_{p+q+2k} - L_{p-q} (-1)^{q+k}  )
\end{equation}
\begin{equation}\label{summand2}
 F_{p+k}F_{q+k}F_{r+k} = \frac{1}{5} ( F_{p+q+r+3k} - (-1)^{r+k} F_{p+q-r+k} - L_{p-q} (-1)^{q+k} F_{r+k}  )
\end{equation}
\begin{equation}\label{summand5}
\begin{split}
 & F_{p+k}F_{q+k}F_{r+k}F_{s+k} \\
 = & \frac{1}{25} ( L_{p+q+r+s+4k} - (-1)^{s+k} L_{p+q+r-s+2k} - (-1)^{r+k} L_{p+q-r+s+2k} \\
 & - L_{p-q} (-1)^{q+k} L_{r+s+2k} + (-1)^{r+s} L_{p+q-r-s} + (-1)^{q+s} L_{p-q}L_{r-s} ) \\
\end{split}
\end{equation}
\begin{equation}\label{summand3}
 L_{p+k}L_{q+k} = L_{p+q+2k} + L_{p-q} (-1)^{q+k}
\end{equation}
\begin{equation}\label{summand4}
 L_{p+k}L_{q+k}L_{r+k} = L_{p+q+r+3k} + (-1)^{r+k} L_{p+q-r+k} + L_{p-q} (-1)^{q+k} L_{r+k}
\end{equation}
\begin{equation}\label{summand6}
\begin{split}
 & L_{p+k}L_{q+k}L_{r+k}L_{s+k} \\
 = &  L_{p+q+r+s+4k} + (-1)^{s+k} L_{p+q+r-s+2k} + (-1)^{r+k} L_{p+q-r+s+2k} \\
 & + L_{p-q} (-1)^{q+k} L_{r+s+2k} + (-1)^{r+s} L_{p+q-r-s} + (-1)^{q+s} L_{p-q}L_{r-s} \\
\end{split}
\end{equation}
When the particular summation identity is also available in this paper,
the summation identity has two evaluations that are also related by the above identities.
For example (\ref{cubic1}) with $G=F$ and $p=q=r=0$ results in:
\begin{equation}
 \sum_{k=0}^n F_k^3 = \frac{1}{2} ( F_{n-1}F_n^2 + F_nF_{n+1}^2 - F_{n+1}F_{n-1}^2 + 1 )
\end{equation}
Applying Cassini's identity $F_{n+1}F_{n-1}-F_n^2=(-1)^n$ this simplifies to:
\begin{equation}
 \sum_{k=0}^n F_k^3 = \frac{1}{2} ( F_nF_{n+1}^2 - (-1)^n F_{n-1} + 1 )
\end{equation}
Using (\ref{summand2}) with $p=q=r=0$:
\begin{equation}\label{fcube}
 F_k^3 = \frac{1}{5} ( F_{3k} - 3(-1)^k F_k )
\end{equation}
and using (10.11) and (10.49) of part 1 of this paper \cite{K19} with $G=F$ and $p=0$:
\begin{equation}
 \sum_{k=0}^n F_{3k} = \frac{1}{2} ( F_{3n+2} - 1 )
\end{equation}
\begin{equation}
 \sum_{k=0}^n (-1)^k F_k = (-1)^n F_{n-1} - 1
\end{equation}
resulting in:
\begin{equation}
 \sum_{k=0}^n F_k^3 = \frac{1}{10} ( F_{3n+2} - 6(-1)^n F_{n-1} + 5 )
\end{equation}
The two right sides of this summation identity are related by (\ref{summand2})
with $p=q=0$ and $r=2$.
Replacing $k$ by $2k$ in (\ref{fcube}) yields:
\begin{equation} 
 F_{2k}^3 = \frac{1}{5} ( F_{6k} - 3 F_{2k} )
\end{equation}
Using (10.7) and (10.17) of part 1 of this paper \cite{K19} with $G=F$ and $p=0$:
\begin{equation}
 \sum_{k=0}^n F_{2k}^3 = \frac{1}{20} ( F_{6n+3} - 12 F_{2n+1} + 10 )
\end{equation}
Using (\ref{summand5}) with $p=q=r=s=0$:
\begin{equation}
 F_k^4 = \frac{1}{25} ( L_{4k} - 4(-1)^k L_{2k} + 6 )
\end{equation}
Using (10.14) and (10.55) of part 1 of this paper \cite{K19} with $G=L$ and $p=0$:
\begin{equation}
 \sum_{k=0}^n F_k^4 = \frac{1}{125} [ L_{4n+4} - L_{4n} - 4(-1)^n ( L_{2n} + L_{2n+2} ) + 30 n + 15 ]
\end{equation}
From (\ref{cubic1}) with $G=L$ and $p=q=r=0$ follows:
\begin{equation}
 \sum_{k=0}^n L_k^3 = \frac{1}{2} ( L_{n-1}L_n^2 + L_nL_{n+1}^2 - L_{n+1}L_{n-1}^2 + 19 )
\end{equation}
For simplifying this identity we need (\ref{cassini}) with $H=G$, $m=n$ and $p=0$:
\begin{equation}
 G_{n+1}G_{n-1} - G_n^2 = (-1)^n (G_1^2-G_0G_2)
\end{equation}
When $G=F$ this results in Cassini's identity. When $G=L$ this results in:
\begin{equation}
 L_{n+1}L_{n-1} - L_n^2 = 5 (-1)^{n+1}
\end{equation}
so the summation identity simplifies to:
\begin{equation}
 \sum_{k=0}^n L_k^3 = \frac{1}{2} ( L_nL_{n+1}^2 + 5(-1)^nL_{n-1} + 19 )
\end{equation}
Likewise using (\ref{summand4}) with $p=q=r=0$:
\begin{equation}
 L_k^3 = L_{3k} + 3(-1)^k L_k
\end{equation}
Using (10.11) and (10.49) of part 1 of this paper \cite{K19} with $G=L$ and $p=0$:
\begin{equation}
 \sum_{k=0}^n L_{3k} = \frac{1}{2} ( L_{3n+2} + 1 )
\end{equation}
\begin{equation}
 \sum_{k=0}^n (-1)^k L_k = (-1)^n L_{n-1} + 3
\end{equation}
resulting in:
\begin{equation}
 \sum_{k=0}^n L_k^3 = \frac{1}{2} ( L_{3n+2} + 6(-1)^n L_{n-1} + 19 )
\end{equation}
The two right sides of this summation identity are related by (\ref{summand4})
with $p=q=0$ and $r=2$.
Using (\ref{summand6}) with $p=q=r=s=0$:
\begin{equation}
 L_k^4 = L_{4k} + 4(-1)^k L_{2k} + 6
\end{equation}
Using (10.14) and (10.55) of part 1 of this paper \cite{K19} with $G=L$ and $p=0$:
\begin{equation}
 \sum_{k=0}^n L_k^4 = \frac{1}{5} [ L_{4n+4} - L_{4n} + 4(-1)^n ( L_{2n} + L_{2n+2} ) + 30 n + 55 ]
\end{equation}
In the list of cubic examples below more examples of this type of summation identities are listed.\\

\section{Derivation of the Quadratic Summation\\ Identities with Binomial Coefficients}

The quadratic summation identities with binomial coefficients can be derived as in the
previous sections, by using again (\ref{prod}) and using two linear summation identities.
The summand of:
\begin{equation}
 \sum_{k=0}^n \binom{n}{k} G_{p+k}G_{q+k}
\end{equation}
is expressed with (\ref{prod}) into:
\begin{equation}
\begin{split}
 G_{p+k}G_{q+k} & = \frac{1}{5} \{ [(3H_0-H_1)G_{p+q+2k} + (2H_1-H_0)G_{p+q+2k+1}] \\
 & + (-1)^q [(2H_0+H_1)G_{p-q} - (2H_1-H_0)G_{p-q+1}] (-1)^k \} \\
\end{split}
\end{equation}
The quadratic sum is expressed in this way as two linear sums, which are evaluated with
part 1 of this paper \cite{K19}, and then the result is evaluated with (\ref{prodsec1}) and (\ref{prodsec2}).
From part 1 of this paper \cite{K19} (11.7):
\begin{equation}
 \sum_{k=0}^n \binom{n}{k} G_{p+q+2k} = 5^{\lfloor n/2\rfloor} ( G_{p+q+n+1} - (-1)^n G_{p+q+n-1} )
\end{equation}
Using \cite{K19} (2.2) with $f(x)=1$ and (3.12), where $\delta_{n,m}$ is the Kronecker delta:
\begin{equation}
\begin{split}
 \sum_{k=0}^n \binom{n}{k} (-1)^k & = \sum_{k=0}^{\infty} (-1)^k \res{x}\frac{(1+x)^n}{x^{k+1}}
 = \res{x} \frac{(1+x)^n}{x} \sum_{k=0}^{\infty} (\frac{-1}{x})^k \\
 & = \res{x} \frac{(1+x)^n}{x(1+1/x)} = \res{x} (1+x)^{n-1} = \delta_{n,0} \\
\end{split}
\end{equation}
Using (\ref{prodsec1}) and (\ref{prodsec2}) the result (\ref{binomres1}) is obtained.
For (\ref{binomres2}) and (\ref{binomres3}) using \cite{K19} (3.13) and (3.14):
\begin{equation}
 \sum_{k=0}^n \binom{n}{k} (-1)^k k = \res{x} \frac{(1+x)^n}{x} \sum_{k=0}^{\infty} k (\frac{-1}{x})^k
 = \res{x} -(1+x)^{n-2} = -\delta_{n,1}
\end{equation}
\begin{equation}
\begin{split}
 \sum_{k=0}^n \binom{n}{k} (-1)^k k^2 & = \res{x} \frac{(1+x)^n}{x} \sum_{k=0}^{\infty} k^2 (\frac{-1}{x})^k \\
 & = \res{x} [ -(1+x)^{n-2} +2(1+x)^{n-3} ] = -\delta_{n,1} + 2 \delta_{n,2} \\
\end{split}
\end{equation}
For identities (12.8) to (12.14) using  \cite{K19} (2.2):
\begin{equation}
\begin{split}
 \sum_{k=0}^n \binom{n}{k} & = \sum_{k=0}^{\infty}  \res{x}\frac{(1+x)^n}{x^{k+1}}
 = \res{x}  \frac{(1+x)^n}{x} \sum_{k=0}^{\infty} (\frac{1}{x})^k \\
 & = \res{x} \frac{(1+x)^n}{x(1-1/x)} = \res{x} \frac{(1+x)^n}{x-1} = 2^n \\
\end{split}
\end{equation}
\begin{equation}
 \sum_{k=0}^n \binom{n}{k} k = \res{x} \frac{(1+x)^n}{x} \sum_{k=0}^{\infty} k (\frac{1}{x})^k
 = \res{x} \frac{(1+x)^n}{(x-1)^2} = n 2^{n-1}
\end{equation}
\begin{equation}
\begin{split}
 \sum_{k=0}^n \binom{n}{k} k^2 & = \res{x} \frac{(1+x)^n}{x} \sum_{k=0}^{\infty} k^2 (\frac{1}{x})^k 
 = \res{x} [ \frac{(1+x)^n}{(x-1)^2} + 2\frac{(1+x)^n}{(x-1)^3} ] \\
 &  = n 2^{n-1} + n(n-1)2^{n-2} = n(n+1)2^{n-2} \\
\end{split}
\end{equation}
\begin{equation}
\end{equation}

\section{Infinite Series and Generating Functions}

When $n\rightarrow\infty$ then from (\ref{gfl}) follows that $G_{a+bn}$ is $O(\phi^{|b|n})$
and $H_{c+dn}$ is $O(\phi^{|d|n})$,
and from (\ref{qres}) follows that $Q_2(v,n)$ is $O(\phi^{(|b|+|d|)n})$.
In that case in the right sides of (\ref{res1}) to (\ref{res4}) therefore $x^n Q_2(v,n)\rightarrow 0$ when:
\begin{equation}
 |x| < \phi^{-|b|-|d|}
\end{equation}
When this condition is fulfilled the following infinite series result
from (\ref{res1}) to (\ref{res6}), with $Q_2(v,w)$ defined in (\ref{qres}).
These are also the generating functions of the product of generalized
Fibonacci numbers $G_{a+bn}H_{c+dn}$.
\begin{equation}\label{infid1}
 \sum_{k=0}^{\infty} x^k G_{a+bk}H_{c+dk} = Q_2(1,0)
\end{equation}
\begin{equation}
 \sum_{k=0}^{\infty} x^k k G_{a+bk}H_{c+dk} = x Q_2(2,1)
\end{equation}
\begin{equation}
 \sum_{k=0}^{\infty} x^k k^2 G_{a+bk}H_{c+dk} = x Q_2(2,1) + 2x^2 Q_2(3,2)
\end{equation}
\begin{equation}
 \sum_{k=0}^{\infty} x^k k^3 G_{a+bk}H_{c+dk} = x Q_2(2,1) + 6x^2 Q_2(3,2) + 6x^3 Q_2(4,3)
\end{equation}
\begin{equation}
 \sum_{k=0}^{\infty} x^k k^4 G_{a+bk}H_{c+dk} = x Q_2(2,1) + 14x^2 Q_2(3,2) + 36x^3 Q_2(4,3) + 24x^4 Q_2(5,4)
\end{equation}
\begin{equation}
\begin{split}
 \sum_{k=0}^{\infty} x^k k^5 G_{a+bk}H_{c+dk} = & x Q_2(2,1) + 30x^2 Q_2(3,2) + 150x^3 Q_2(4,3) \\
 & + 240x^4 Q_2(5,4) + 120x^5 Q_2(6,5) \\
\end{split}
\end{equation}
From the results of section 2,
the following are quadratic infinite series for $|x|<\phi^{-2}$, which are also generating functions
of the product of generalized Fibonacci numbers $G_{p+n}H_{q+n}$ and $G_{p+n}H_{q-n}$.
\begin{equation}
\begin{split}
 \sum_{k=0}^{\infty} & x^k G_{p+k}H_{q+k} = \frac{1}{1-3x+x^2} ( G_pH_q - x G_{p-1}H_{q-1} ) \\
 & - \frac{x}{(1+x)(1-3x+x^2)} ( 3G_pH_q - G_{p-1}H_{q-1} - G_{p+1}H_{q+1} ) \\
\end{split}
\end{equation}
\begin{equation}
\begin{split}
 \sum_{k=0}^{\infty} & x^k G_{p+k}H_{q-k} = \frac{1}{1+3x+x^2} ( G_pH_q - x G_{p-1}H_{q+1} ) \\
 & + \frac{x}{(1-x)(1+3x+x^2)} ( 2G_pH_q + G_{p-1}H_{q-1} + G_{p+1}H_{q+1} ) \\
\end{split}
\end{equation}
From the results of section 3,
the following are cubic infinite series for $|x|<\phi^{-3}$, which are also generating functions
of the product of generalized Fibonacci numbers $G_{p+n}H_{q+n}K_{r+n}$,  $G_{p+n}H_{q+n}K_{r-n}$
 and  $G_{p+n}H_{q-n}K_{r-n}$.
\begin{equation}\label{gencube1}
\begin{split}
 \sum_{k=0}^{\infty} x^k G_{p+k}H_{q+k}K_{r+k} = & \frac{-1}{(1+x-x^2)(1-4x-x^2)} \\ 
 \cdot [ & x^2(x+3) G_{p-1}H_{q-1}K_{r-1} + (3x-1) G_pH_qK_r \\
 & + x^2 G_{p-2}H_{q-2}K_{r-2} - x G_{p+1}H_{q+1}K_{r+1} ] \\
\end{split}
\end{equation}
\begin{equation}\label{gencube2}
\begin{split}
 \sum_{k=0}^{\infty} x^k G_{p+k}H_{q+k}K_{r-k} = & \frac{-1}{(1-x-x^2)(1+4x-x^2)} \\ 
 \cdot [ & x^2(x-3) G_{p-1}H_{q-1}K_{r+1} - (3x+1) G_pH_qK_r \\
 & + x^2 G_{p-2}H_{q-2}K_{r+2} - x G_{p+1}H_{q+1}K_{r-1} ] \\
\end{split}
\end{equation}
\begin{equation}\label{gencube3}
\begin{split}
 \sum_{k=0}^{\infty} x^k G_{p+k}H_{q-k}K_{r-k} = & \frac{-1}{(1+x-x^2)(1-4x-x^2)} \\ 
 \cdot [ & x^2(x+3) G_{p-1}H_{q+1}K_{r+1} + (3x-1) G_pH_qK_r \\
 & + x^2 G_{p-2}H_{q+2}K_{r+2} - x G_{p+1}H_{q-1}K_{r-1} ] \\
\end{split}
\end{equation}
and the following for $|x|<\phi^{-6}$:
\begin{equation}\label{gencube4}
\begin{split}
 \sum_{k=0}^{\infty} x^k G_{p+2k}H_{q+2k}K_{r+2k} = & \frac{-1}{(1-3x+x^2)(1-18x+x^2)} \\ 
 \cdot [ & x^2(x-21) G_{p-2}H_{q-2}K_{r-2} + (21x-1) G_pH_qK_r \\
 & + x^2 G_{p-4}H_{q-4}K_{r-4} - x G_{p+2}H_{q+2}K_{r+2} ] \\
\end{split}
\end{equation}
Taking $G=H=K=F$ or $G=H=K=L$ these formulas include the quadratic and cubic generating functions listed in \cite{HL}.

\section{Simplification of Sums of Products of\\ Generalized Fibonacci Numbers}

For the simplification of sums of products of generalized Fibonacci numbers,
the following equation is solved by a computer program:
\begin{equation}
 \sum_{k=0}^n a_k \prod_{j=1}^p F_{q_j+c_{k,j}} = \sum_{k=0}^m b_k \prod_{j=1}^p F_{q_j+d_{k,j}}
\end{equation}
In this formula the left side is the  expression that needs to be simplified,
where $m<n$ and the parameters $a_k$ and $c_{k,j}$ are constant input parameters,
and the $q_j$ are any integers but in the program they can be restricted to a few,
for example $-2\leq q_j\leq 2$. The $b_k$ and $d_{k,j}$ are generated
within some range and the identity is checked for all values of $q_j$.
For efficiency the result of the left side is stored for all values of $q_j$,
and for the right side when the $d_{k,j}$ are generated the products
are stored for all values of $q_j$ before the $b_k$ are generated.
When the identity is true for all values of $q_j$ the parameters
$b_k$ and $d_{k,j}$ are listed.\\
For this computer program the Fibonacci number $F_n$ must be computed.
For this the following identities are used:
\begin{equation}
 F_{n+1} = \frac{1}{2} ( F_n + L_n )
\end{equation}
\begin{equation}
 L_{n+1} = \frac{1}{2} ( 5 F_n + L_n )
\end{equation}
\begin{equation}
 F_{2n} = F_n L_n
\end{equation}
\begin{equation}
 L_{2n} = L_n^2 - 2 (-1)^n
\end{equation}
Starting with $F_0=0$ and $L_0=2$,
for every bit of $n$ in decreasing order the first and second identities are computed when
the bit is one, and then the third and fourth identities when there are bits left,
thus requiring $\lfloor\log_2(n)\rfloor+1$ iterations.
The complexity of this algorithm is therefore order $\log(n)M(n)$ where $M(n)$ is the
complexity of multiplication.
When $n$ is negative the identity $F_{-n}=(-1)^{n+1}F_n$ is used.
For computing generalized Fibonacci numbers the same algorithm is used
with (\ref{gfl}) and (\ref{gmin}).

\section{List of Quadratic Examples}

Examples with $x=1$, $b>0$ and $d>0$:
\begin{equation}\label{quad1}
 \sum_{k=0}^n G_{p+k}H_{q+k} = \frac{1}{2} ( G_{p+n}H_{q+n+1} + G_{p+n+1}H_{q+n} - G_{p-1}H_q - G_p H_{q-1} )
\end{equation}
\begin{equation}
\begin{split}
 \sum_{k=0}^n & k G_{p+k}H_{q+k} = \frac{1}{2} n ( G_{p+n}H_{q+n+1} + G_{p+n+1}H_{q+n} ) \\
 & - \frac{1}{4} ( G_{p+n+1}H_{q+n+2} - G_{p+n-1}H_{q+n-2} - G_{p+1}H_{q+2} + G_{p-1}H_{q-2} ) \\
\end{split}
\end{equation}
\begin{equation}
\begin{split}
 \sum_{k=0}^n k^2 G_{p+k}H_{q+k} = & \frac{1}{2} n^2 ( G_{p+n}H_{q+n+1} + G_{p+n+1}H_{q+n} ) \\
 & - \frac{1}{2} n ( G_{p+n+1}H_{q+n+2} - G_{p+n-1}H_{q+n-2} ) \\
 & + G_{p+n}H_{q+n-1} + G_{p+n+1}H_{q+n} - G_pH_{q-1} - G_{p+1}H_q  \\
\end{split}
\end{equation}
\begin{equation}
\begin{split}
 \sum_{k=0}^n & k^3 G_{p+k}H_{q+k} = \frac{1}{2} n^3 ( G_{p+n}H_{q+n+1} + G_{p+n+1}H_{q+n} ) \\
 & - \frac{3}{4} n^2 ( G_{p+n+1}H_{q+n+2} - G_{p+n-1}H_{q+n-2} ) \\
 & + 3 n ( G_{p+n}H_{q+n-1} + G_{p+n+1}H_{q+n} ) - \frac{1}{8} [ 23 ( G_{p+n}H_{q+n} - G_pH_q ) \\
 & + 11( G_{p+n-1}H_{q+n-1} + G_{p+n+1}H_{q+n+1} - G_{p-1}H_{q-1} - G_{p+1}H_{q+1} ) ]  \\
\end{split}
\end{equation}
\begin{equation}
 \sum_{k=0}^n G_{p+2k}H_{q+k} = \frac{1}{2} ( G_{p+2n+1}H_{q+n} + G_{p+2n-1}H_{q+n+1} + G_p H_q - G_{p-1}H_{q+2} )
\end{equation}
\begin{equation}
\begin{split}
 \sum_{k=0}^n & k G_{p+2k}H_{q+k} = \frac{1}{2} n ( G_{p+2n+1}H_{q+n} + G_{p+2n-1}H_{q+n+1} ) \\
 & - \frac{1}{4} [ 3 ( G_{p+2n-2}H_{q+n+1} - G_{p-2}H_{q+1} ) - G_{p+2n-3}H_{q+n} + G_{p-3}H_q ] \\
\end{split}
\end{equation}
\begin{equation}
\begin{split}
 & \sum_{k=0}^n k^2 G_{p+2k}H_{q+k} = \frac{1}{2} n^2 ( G_{p+2n+1}H_{q+n} + G_{p+2n-1}H_{q+n+1} ) \\
 & - \frac{1}{2} n ( 3G_{p+2n-2}H_{q+n+1} - G_{p+2n-3}H_{q+n} ) \\
 & + \frac{1}{4} [ 3 ( G_{p+2n-3}H_{q+n+3} - G_{p-3}H_{q+3} ) - G_{p+2n-4}H_{q+n+2}  + G_{p-4}H_{q+2} ] \\
\end{split}
\end{equation}
\begin{equation}
\begin{split}
 \sum_{k=0}^n & G_{p+2k}H_{q+2k} = \frac{1}{5} [ G_{p+2n}H_{q+2n+1} + G_{p+2n+1}H_{q+2n+2} \\
  & + n ( 3G_pH_q - G_{p-1}H_{q-1} - G_{p+1}H_{q+1} ) + 6G_pH_q - G_{p+2}H_{q+2} ] \\
\end{split}
\end{equation}
\begin{equation}
\begin{split}
 \sum_{k=0}^n & k G_{p+2k}H_{q+2k} = \frac{1}{5} n ( G_{p+2n}H_{q+2n+1} + G_{p+2n+1}H_{q+2n+2} ) \\
  & - \frac{1}{25} ( 2G_{p+2n}H_{q+2n} + G_{p+2n-1}H_{q+2n-1} + G_{p+2n+1}H_{q+2n+1} \\
  & \qquad - 2G_pH_q - G_{p-1}H_{q-1} - G_{p+1}H_{q+1} ) \\
  & + \frac{1}{10} n(n+1) ( 3G_pH_q - G_{p-1}H_{q-1} - G_{p+1}H_{q+1} ) \\
\end{split}
\end{equation}
\begin{equation}
\begin{split}
 \sum_{k=0}^n & k^2 G_{p+2k}H_{q+2k} = \frac{1}{5} n^2 ( G_{p+2n}H_{q+2n+1} + G_{p+2n+1}H_{q+2n+2} ) \\
  & - \frac{1}{25} [ 2n ( 2G_{p+2n}H_{q+2n} + G_{p+2n-1}H_{q+2n-1} + G_{p+2n+1}H_{q+2n+1} ) \\
  & \quad - 3 ( G_{p+2n-1}H_{q+2n} + G_{p+2n}H_{q+2n+1} - G_{p-1}H_q - G_pH_{q+1} ) ] \\
  & + \frac{1}{30} n(n+1)(2n+1) ( 3G_pH_q - G_{p-1}H_{q-1} - G_{p+1}H_{q+1} ) \\
\end{split}
\end{equation}
\begin{equation}
 \sum_{k=0}^n G_{p+3k}H_{q+k} = \frac{1}{5} ( G_{p+3n}H_{q+n+1} + G_{p+3n+3}H_{q+n} - G_pH_{q-1} - G_{p-3}H_q )
\end{equation}
\begin{equation}
 \sum_{k=0}^n G_{p+3k}H_{q+3k} = \frac{1}{4} ( G_{p+3n}H_{q+3n+1} + G_{p+3n+2}H_{q+3n}
   + G_pH_{q-2} - G_{p-1}H_q ) 
\end{equation}
\begin{equation}
\begin{split}
 \sum_{k=0}^n G_{p+4k}H_{q+4k} = & \frac{1}{15} [ G_{p+4n+2}H_{q+4n+1} + G_{p+4n+3}H_{q+4n+2} \\
  & + 3n ( 3G_pH_q - G_{p-1}H_{q-1} - G_{p+1}H_{q+1} ) \\
  & + 12G_pH_q - 5G_{p+1}H_{q+1} + 2G_{p-1}H_{q-1} ] \\
\end{split}
\end{equation}
\begin{equation}
 \sum_{k=0}^n G_{p+5k}H_{q+5k} = \frac{1}{22} ( G_{p+5n}H_{q+5n+5} + G_{p+5n+5}H_{q+5n}
   - G_pH_{q-5} - G_{p-5}H_q ) 
\end{equation}
Examples with $x=1$, $b>0$ and $d<0$:
\begin{equation}\label{resultconv1}
\begin{split}
 \sum_{k=0}^n & G_{p+k}H_{q-k} = \frac{1}{5} [ G_{p+n}H_{q-n} - G_{p+n+1}H_{q-n-1}  \\
  & + n ( 2G_pH_q + G_{p-1}H_{q-1} + G_{p+1}H_{q+1} ) + 4 G_pH_q + G_{p+1}H_{q-1} ] \\
\end{split}
\end{equation}
\begin{equation}
\begin{split}
 \sum_{k=0}^n & k G_{p+k}H_{q-k} = \frac{1}{5} n ( G_{p+n}H_{q-n} - G_{p+n+1}H_{q-n-1} ) \\
  & + \frac{1}{25} ( 3G_{p+n}H_{q-n} - G_{p+n-1}H_{q-n-1} - G_{p+n+1}H_{q-n+1}  \\
  & \qquad - 3G_pH_q + G_{p-1}H_{q-1} + G_{p+1}H_{q+1} ) \\
  & + \frac{1}{10} n(n+1) ( 2G_pH_q + G_{p-1}H_{q-1} + G_{p+1}H_{q+1} ) \\
\end{split}
\end{equation}
\begin{equation}
\begin{split}
 \sum_{k=0}^n & k^2 G_{p+k}H_{q-k} = \frac{1}{5} n^2 ( G_{p+n}H_{q-n} - G_{p+n+1}H_{q-n-1} ) \\
  & + \frac{1}{25} [ 2n ( 3G_{p+n}H_{q-n} - G_{p+n-1}H_{q-n-1} - G_{p+n+1}H_{q-n+1} )  \\
  & \qquad + G_{p+n}H_{q-n-1} - G_{p+n-1}H_{q-n} - G_pH_{q-1} + G_{p-1}H_q ] \\
  & + \frac{1}{30} n(n+1)(2n+1) ( 2G_pH_q + G_{p-1}H_{q-1} + G_{p+1}H_{q+1} ) \\
\end{split}
\end{equation}
\begin{equation}
\begin{split}
 \sum_{k=0}^n & k^3 G_{p+k}H_{q-k} = \frac{1}{5} n^3 ( G_{p+n}H_{q-n} - G_{p+n+1}H_{q-n-1} ) \\
  & + \frac{3}{25} n [ n ( 3G_{p+n}H_{q-n} - G_{p+n-1}H_{q-n-1} - G_{p+n+1}H_{q-n+1} )  \\
  & \qquad + G_{p+n}H_{q-n-1} - G_{p+n-1}H_{q-n}  ] \\
  & - \frac{1}{125} ( 3G_{p+n}H_{q-n} - G_{p+n-1}H_{q-n-1} - G_{p+n+1}H_{q-n+1} \\
  & \qquad - 3G_pH_q + G_{p-1}H_{q-1} + G_{p+1}H_{q+1} ) \\
  & + \frac{1}{20} n^2(n+1)^2 ( 2G_pH_q + G_{p-1}H_{q-1} + G_{p+1}H_{q+1} ) \\
\end{split}
\end{equation}
\begin{equation}
 \sum_{k=0}^n G_{p+2k}H_{q-k} =  \frac{1}{2} ( G_{p+2n+2}H_{q-n} + G_{p+2n+1}H_{q-n+1} 
  - G_{p-1}H_q - G_{p+1}H_{q+1} )  
\end{equation}
\begin{equation}
 \sum_{k=0}^n G_{p+k}H_{q-2k} = \frac{1}{2} ( G_{p+n}H_{q-2n} - G_{p+n+2}H_{q-2n+1}
  + G_pH_q + G_{p+2}H_{q+1} ) 
\end{equation}
\begin{equation}
\begin{split}
 \sum_{k=0}^n & G_{p+2k}H_{q-2k} = \frac{1}{5} [ G_{p+2n+1}H_{q-2n} - G_{p+2n+2}H_{q-2n-1}  \\
  & + n ( 2G_pH_q + G_{p-1}H_{q-1} + G_{p+1}H_{q+1} ) + 6 G_pH_q - G_{p+2}H_{q-2} ] \\
\end{split}
\end{equation}
\begin{equation}
\begin{split}
 \sum_{k=0}^n & k G_{p+2k}H_{q-2k} = \frac{1}{5} n ( G_{p+2n+1}H_{q-2n} - G_{p+2n+2}H_{q-2n-1} )  \\
  & - \frac{1}{25} ( 3G_{p+2n}H_{q-2n} - G_{p+2n-1}H_{q-2n-1} - G_{p+2n+1}H_{q-2n+1} \\
  & \qquad - 3G_pH_q + G_{p-1}H_{q-1} + G_{p+1}H_{q+1} ) \\
  & + \frac{1}{10} n(n+1) ( 2G_pH_q + G_{p-1}H_{q-1} + G_{p+1}H_{q+1} ) \\
\end{split}
\end{equation}
\begin{equation}
\begin{split}
 \sum_{k=0}^n & k^2 G_{p+2k}H_{q-2k} = \frac{1}{5} n^2 ( G_{p+2n+1}H_{q-2n} - G_{p+2n+2}H_{q-2n-1} )  \\
  & - \frac{1}{25} [ 2n ( 3G_{p+2n}H_{q-2n} - G_{p+2n-1}H_{q-2n-1} - G_{p+2n+1}H_{q-2n+1} ) \\
  & \qquad - 3 ( G_{p+2n-1}H_{q-2n} - G_{p+2n}H_{q-2n-1} - G_{p-1}H_q + G_pH_{q-1} ) ] \\
  & + \frac{1}{30} n(n+1)(2n+1) ( 2G_pH_q + G_{p-1}H_{q-1} + G_{p+1}H_{q+1} ) \\
\end{split}
\end{equation}
\begin{equation}
 \sum_{k=0}^n G_{p+3k}H_{q-k} =  \frac{1}{3} ( G_{p+3n+3}H_{q-n} + G_{p+3n}H_{q-n-1} 
  - G_{p-3}H_q - G_pH_{q+1} )  
\end{equation}
\begin{equation}
 \sum_{k=0}^n G_{p+k}H_{q-3k} =  - \frac{1}{3} ( G_{p+n+1}H_{q-3n} + G_{p+n}H_{q-3n-3} 
  - G_pH_{q+3} - G_{p-1}H_q )  
\end{equation}
\begin{equation}
\begin{split}
 \sum_{k=0}^n G_{p+3k}H_{q-3k} = & \frac{1}{10} [ G_{p+3n+2}H_{q-3n} - G_{p+3n+3}H_{q-3n-1} \\
  & + 2n ( 2G_pH_q + G_{p-1}H_{q-1} + G_{p+1}H_{q+1} ) \\
  & + 8G_pH_q + 3G_{p+1}H_{q-1} - G_{p-1}H_{q+1} ] \\
\end{split}
\end{equation}
\begin{equation}
\begin{split}
 \sum_{k=0}^n G_{p+4k}H_{q-4k} = & \frac{1}{15} [ G_{p+4n+3}H_{q-4n} - G_{p+4n+4}H_{q-4n-1} \\
  & + 3n ( 2G_pH_q + G_{p-1}H_{q-1} + G_{p+1}H_{q+1} ) \\
  & + 12G_pH_q + 5G_{p+1}H_{q-1} - 2G_{p-1}H_{q+1} ] \\
\end{split}
\end{equation}
\begin{equation}
\begin{split}
 \sum_{k=0}^n G_{p+5k}H_{q-5k} = & \frac{1}{25} [ G_{p+5n+4}H_{q-5n} - G_{p+5n+5}H_{q-5n-1} \\
  & + 5n ( 2G_pH_q + G_{p-1}H_{q-1} + G_{p+1}H_{q+1} ) \\
  & + 20G_pH_q + 8G_{p+1}H_{q-1} - 3G_{p-1}H_{q+1} ] \\
\end{split}
\end{equation}
Examples with $x=-1$, $b>0$ and $d>0$:
\begin{equation}\label{reslist1}
\begin{split}
 \sum_{k=0}^n & (-1)^k G_{p+k}H_{q+k} = \frac{1}{5} [ (-1)^n ( G_{p+n}H_{q+n} + G_{p+n+1}H_{q+n+1} ) \\
  & + n ( 3 G_p H_q - G_{p-1}H_{q-1} - G_{p+1}H_{q+1} ) + 4 G_p H_q - G_{p+1}H_{q+1}  ] \\
\end{split}
\end{equation}
\begin{equation}
\begin{split}
 \sum_{k=0}^n & (-1)^k k G_{p+k}H_{q+k} = \frac{1}{5} (-1)^n n ( G_{p+n}H_{q+n} + G_{p+n+1}H_{q+n+1} ) \\
  & + \frac{1}{25} [ (-1)^n ( 2G_{p+n}H_{q+n} + G_{p+n-1}H_{q+n-1} + G_{p+n+1}H_{q+n+1} ) \\
  & \quad - 2G_pH_q - G_{p-1}H_{q-1} - G_{p+1}H_{q+1} ] \\
  & + \frac{1}{10} n(n+1) ( 3G_pH_q - G_{p-1}H_{q-1} - G_{p+1}H_{q+1} ) \\
\end{split}
\end{equation}
\begin{equation}
\begin{split}
 \sum_{k=0}^n & (-1)^k k^2 G_{p+k}H_{q+k} = \frac{1}{5} (-1)^n n^2 ( G_{p+n}H_{q+n} + G_{p+n+1}H_{q+n+1} ) \\
  & + \frac{1}{25} \{ (-1)^n [ 2n ( 2G_{p+n}H_{q+n} + G_{p+n-1}H_{q+n-1} + G_{p+n+1}H_{q+n+1} ) \\
  & \quad - G_{p+n}H_{q+n-1} - G_{p+n+1}H_{q+n} ] + G_pH_{q-1} + G_{p+1}H_q \} \\
  & + \frac{1}{30} n(n+1)(2n+1) ( 3G_pH_q - G_{p-1}H_{q-1} - G_{p+1}H_{q+1} ) \\
\end{split}
\end{equation}
\begin{equation}
\begin{split}
 \sum_{k=0}^n & (-1)^k k^3 G_{p+k}H_{q+k} = \frac{1}{5} (-1)^n n^3 ( G_{p+n}H_{q+n} + G_{p+n+1}H_{q+n+1} ) \\
  & + \frac{3}{25} (-1)^n n [ n ( 2G_{p+n}H_{q+n} + G_{p+n-1}H_{q+n-1} + G_{p+n+1}H_{q+n+1} ) \\
  & \qquad - G_{p+n}H_{q+n-1} - G_{p+n+1}H_{q+n} ] \\
  & - \frac{1}{125} [ (-1)^n ( 2G_{p+n}H_{q+n} + G_{p+n-1}H_{q+n-1} + G_{p+n+1}H_{q+n+1} ) \\
  & \qquad - 2G_pH_q - G_{p-1}H_{q-1} - G_{p+1}H_{q+1} ] \\
  & + \frac{1}{20} n^2(n+1)^2 ( 3G_pH_q - G_{p-1}H_{q-1} - G_{p+1}H_{q+1} ) \\
\end{split}
\end{equation}
\begin{equation}
\begin{split}
 \sum_{k=0}^n (-1)^k G_{p+2k}H_{q+k} = \frac{1}{2} & [ (-1)^n ( G_{p+2n+3}H_{q+n} - G_{p+2n+1}H_{q+n+1} ) \\
  & + G_{p+1}H_{q-1} - G_{p-1}H_q ] \\
\end{split}
\end{equation}
\begin{equation}
\begin{split}
 \sum_{k=0}^n & (-1)^k k G_{p+2k}H_{q+k} = \frac{1}{2} (-1)^n n ( G_{p+2n+3}H_{q+n} - G_{p+2n+1}H_{q+n+1} ) \\
 & + \frac{1}{4} [ (-1)^n ( 3G_{p+2n+2}H_{q+n-1} - G_{p+2n+3}H_{q+n} ) - 3G_{p+2}H_{q-1} + G_{p+3}H_q ] \\
\end{split}
\end{equation}
\begin{equation}
\begin{split}
 \sum_{k=0}^n (-1)^k k^2 G_{p+2k}H_{q+k} = & \frac{1}{2} (-1)^n n [ n ( G_{p+2n+3}H_{q+n} - G_{p+2n+1}H_{q+n+1} ) \\
 & \quad + 3G_{p+2n+2}H_{q+n-1} - G_{p+2n+3}H_{q+n} ] \\
 & - \frac{1}{4} [ (-1)^n ( 3G_{p+2n+3}H_{q+n-3} - G_{p+2n+4}H_{q+n-2} ) \\
 & \quad - 3G_{p+3}H_{q-3} + G_{p+4}H_{q-2} ] \\
\end{split}
\end{equation}
\begin{equation}
\begin{split}
 \sum_{k=0}^n (-1)^k G_{p+2k}H_{q+2k} = & \frac{1}{6} [ (-1)^n ( G_{p+2n}H_{q+2n+2} + G_{p+2n+2}H_{q+2n} ) \\
 & + G_pH_{q-2} + G_{p-2}H_q ] \\
\end{split}
\end{equation}
\begin{equation}
\begin{split}
 \sum_{k=0}^n & (-1)^k k G_{p+2k}H_{q+2k} = \frac{1}{6} (-1)^n n ( G_{p+2n}H_{q+2n+2} + G_{p+2n+2}H_{q+2n} ) \\
 & + \frac{1}{36} [ (-1)^n ( 7G_{p+2n}H_{q+2n} - G_{p+2n-1}H_{q+2n-1} - G_{p+2n+1}H_{q+2n+1} ) \\
 & \quad - 7G_pH_q + G_{p-1}H_{q-1} + G_{p+1}H_{q+1} ] \\
\end{split}
\end{equation}
\begin{equation}
\begin{split}
 \sum_{k=0}^n & (-1)^k k^2 G_{p+2k}H_{q+2k} = \frac{1}{6} (-1)^n n^2 ( G_{p+2n}H_{q+2n+2} + G_{p+2n+2}H_{q+2n} ) \\
 & + \frac{1}{18} (-1)^n n ( 7G_{p+2n}H_{q+2n} - G_{p+2n-1}H_{q+2n-1} - G_{p+2n+1}H_{q+2n+1} ) \\
 & - \frac{1}{27} [ (-1)^n ( G_{p+2n}H_{q+2n-1} + G_{p+2n+1}H_{q+2n} ) - G_pH_{q-1} - G_{p+1}H_q ] \\
\end{split}
\end{equation}
\begin{equation}
\begin{split}
 \sum_{k=0}^n (-1)^k G_{p+3k}H_{q+k} = & \frac{1}{3} [ (-1)^n ( G_{p+3n+3}H_{q+n} - G_{p+3n}H_{q+n+1} ) \\
 & \quad + G_pH_{q-1} - G_{p-3}H_q ] \\
\end{split}
\end{equation}
\begin{equation}
\begin{split}
 \sum_{k=0}^n (-1)^k G_{p+3k}H_{q+3k} = & \frac{1}{10} (-1)^n ( G_{p+3n+1}H_{q+3n+1} + G_{p+3n+2}H_{q+3n+2} ) \\
 & + \frac{1}{5} n ( 3G_pH_q - G_{p-1}H_{q-1} - G_{p+1}H_{q+1} ) \\
 & + \frac{1}{10} ( 5G_pH_q - 3G_pH_{q-1} - G_{p-1}H_{q+3} ) \\
\end{split}
\end{equation}
\begin{equation}
\begin{split}
 \sum_{k=0}^n (-1)^k G_{p+4k}H_{q+4k} = & \frac{1}{14} [ (-1)^n ( G_{p+4n}H_{q+4n+4} + G_{p+4n+4}H_{q+4n} ) \\
 & \quad + G_pH_{q-4} + G_{p-4}H_q ] \\
\end{split}
\end{equation}
\begin{equation}
\begin{split}
 \sum_{k=0}^n (-1)^k G_{p+5k}H_{q+5k} = & \frac{1}{25} (-1)^n ( G_{p+5n+2}H_{q+5n+2} + G_{p+5n+3}H_{q+5n+3} ) \\
 & + \frac{1}{5} n ( 3G_pH_q - G_{p-1}H_{q-1} - G_{p+1}H_{q+1} ) \\
 & + \frac{1}{25} ( 12G_pH_q - 8G_pH_{q-1} - G_{p-1}H_{q+5} ) \\
\end{split}
\end{equation}
Examples with $x=-1$, $b>0$ and $d<0$:
\begin{equation}
\begin{split}
 \sum_{k=0}^n (-1)^k G_{p+k}H_{q-k} = \frac{1}{2} & [ (-1)^n ( G_{p+n+1}H_{q-n} - G_{p+n}H_{q-n-1} ) \\
 &  - G_{p-1}H_q + G_p H_{q+1} ] \\
\end{split}
\end{equation}
\begin{equation}
\begin{split}
 \sum_{k=0}^n & (-1)^k k G_{p+k}H_{q-k} = \frac{1}{2} (-1)^n n ( G_{p+n+1}H_{q-n} - G_{p+n}H_{q-n-1} ) \\
 & - \frac{1}{4} [ (-1)^n ( 2G_{p+n}H_{q-n} - G_{p+n-1}H_{q-n-1} - G_{p+n+1}H_{q-n+1} ) \\
 & \quad - 2G_pH_q + G_{p-1}H_{q-1} + G_{p+1}H_{q+1} ] \\
\end{split}
\end{equation}
\begin{equation}
\begin{split}
 \sum_{k=0}^n & (-1)^k k^2 G_{p+k}H_{q-k} = \frac{1}{2} (-1)^n n [ n ( G_{p+n+1}H_{q-n} - G_{p+n}H_{q-n-1} ) \\
 & \quad - 2G_{p+n}H_{q-n} + G_{p+n-1}H_{q-n-1} + G_{p+n+1}H_{q-n+1} ] \\
 & + (-1)^n ( G_{p+n-1}H_{q-n} - G_{p+n}H_{q-n-1} ) - G_{p-1}H_q + G_pH_{q-1} \\
\end{split}
\end{equation}
\begin{equation}
\begin{split}
 \sum_{k=0}^n & (-1)^k k^3 G_{p+k}H_{q-k} = \frac{1}{2} (-1)^n n^3 ( G_{p+n+1}H_{q-n} - G_{p+n}H_{q-n-1} ) \\
 & - \frac{3}{4} (-1)^n n^2 ( 2G_{p+n}H_{q-n} - G_{p+n-1}H_{q-n-1} - G_{p+n+1}H_{q-n+1} ) \\
 & + 3 (-1)^n n ( G_{p+n-1}H_{q-n} - G_{p+n}H_{q-n-1} ) \\
 & - \frac{1}{8} [ (-1)^n ( 34G_{p+n}H_{q-n} - 11G_{p+n-1}H_{q-n-1} - 11G_{p+n+1}H_{q-n+1} ) \\
 & \quad - 34G_pH_q + 11G_{p-1}H_{q-1} + 11G_{p+1}H_{q+1} ] \\
\end{split}
\end{equation}
\begin{equation}
\begin{split}
 \sum_{k=0}^n (-1)^k G_{p+2k}H_{q-k} = \frac{1}{2} & [ (-1)^n ( G_{p+2n+1}H_{q-n} - G_{p+2n-1}H_{q-n-1} ) \\
 &  + G_{p-2}H_q + G_{p-1}H_{q-1} ] \\
\end{split}
\end{equation}
\begin{equation}
\begin{split}
 \sum_{k=0}^n (-1)^k G_{p+k}H_{q-2k} = \frac{1}{2} & [ (-1)^n ( G_{p+n}H_{q-2n} - G_{p+n-2}H_{q-2n-1} ) \\
 &  + G_pH_q + G_{p-2}H_{q-1} ] \\
\end{split}
\end{equation}
\begin{equation}
\begin{split}
 \sum_{k=0}^n (-1)^k G_{p+2k}H_{q-2k} = \frac{1}{6} & [ (-1)^n ( G_{p+2n+2}H_{q-2n} + G_{p+2n}H_{q-2n-2} ) \\
 &  + G_pH_{q+2} + G_{p-2}H_q ] \\
\end{split}
\end{equation}
\begin{equation}
\begin{split}
 \sum_{k=0}^n & (-1)^k k G_{p+2k}H_{q-2k} = \frac{1}{6} (-1)^n n ( G_{p+2n+2}H_{q-2n} + G_{p+2n}H_{q-2n-2} ) \\
 & + \frac{1}{36} [ (-1)^n ( 6G_{p+2n}H_{q-2n} + G_{p+2n-1}H_{q-2n-1} + G_{p+2n+1}H_{q-2n+1} ) \\
 & \quad - 6G_pH_q - G_{p-1}H_{q-1} - G_{p+1}H_{q+1} ] \\
\end{split}
\end{equation}
\begin{equation}
\begin{split}
 \sum_{k=0}^n & (-1)^k k^2 G_{p+2k}H_{q-2k} = \frac{1}{6} (-1)^n n^2 ( G_{p+2n+2}H_{q-2n} + G_{p+2n}H_{q-2n-2} ) \\
 & + \frac{1}{18} (-1)^n n ( 6G_{p+2n}H_{q-2n} + G_{p+2n-1}H_{q-2n-1} + G_{p+2n+1}H_{q-2n+1} ) \\
 & - \frac{1}{27} [ (-1)^n ( G_{p+2n-1}H_{q-2n} - G_{p+2n}H_{q-2n-1} ) - G_{p-1}H_q + G_pH_{q-1} ] \\
\end{split}
\end{equation}
\begin{equation}
\begin{split}
 \sum_{k=0}^n (-1)^k G_{p+3k}H_{q-k} = \frac{1}{5} & [ (-1)^n ( G_{p+3n+3}H_{q-n} - G_{p+3n}H_{q-n-1} ) \\
 &  + G_pH_{q+1} - G_{p-3}H_q ] \\
\end{split}
\end{equation}
\begin{equation}
\begin{split}
 \sum_{k=0}^n (-1)^k G_{p+k}H_{q-3k} = \frac{1}{5} & [ (-1)^n ( G_{p+n+1}H_{q-3n} - G_{p+n}H_{q-3n-3} ) \\
 &  + G_pH_{q+3} - G_{p-1}H_q ] \\
\end{split}
\end{equation}
\begin{equation}
\begin{split}
 \sum_{k=0}^n (-1)^k G_{p+3k}H_{q-3k} = \frac{1}{4} & [ (-1)^n ( G_{p+3n+1}H_{q-3n} + G_{p+3n}H_{q-3n-2} ) \\
 &  + G_pH_{q+1} + G_{p-2}H_q ] \\
\end{split}
\end{equation}
\begin{equation}
\begin{split}
 \sum_{k=0}^n (-1)^k G_{p+4k}H_{q-4k} = \frac{1}{14} & [ (-1)^n ( G_{p+4n+4}H_{q-4n} + G_{p+4n}H_{q-4n-4} ) \\
 &  + G_pH_{q+4} + G_{p-4}H_q ] \\
\end{split}
\end{equation}
\begin{equation}
\begin{split}
 \sum_{k=0}^n (-1)^k G_{p+5k}H_{q-5k} = \frac{1}{22} & [ (-1)^n ( G_{p+5n+5}H_{q-5n} - G_{p+5n}H_{q-5n-5} ) \\
 &  + G_pH_{q+5} - G_{p-5}H_q ] \\
\end{split}
\end{equation}
Examples with $x=2$:
\begin{equation}
\begin{split}
 \sum_{k=0}^n 2^k G_{p+k}H_{q+k} = \frac{1}{3} & [ 2^{n+1} ( G_{p+n}H_{q+n+1} + G_{p+n-1}H_{q+n-2} ) \\
 & - G_{p-1}H_q - G_{p-2}H_{q-3} ] \\
\end{split}
\end{equation}
\begin{equation}
\begin{split}
 \sum_{k=0}^n & 2^k k G_{p+k}H_{q+k} = \frac{1}{3} 2^{n+1} n ( G_{p+n}H_{q+n+1} + G_{p+n-1}H_{q+n-2} ) \\
 & + \frac{2}{9} [ 2^n ( 2G_{p+n}H_{q+n} - 6G_{p+n-1}H_{q+n-1} - 7G_{p+n-2}H_{q+n-2} ) \\
 & \quad - 2G_pH_q + 6G_{p-1}H_{q-1} + 7G_{p-2}H_{q-2} ] \\
\end{split}
\end{equation}
\begin{equation}
\begin{split}
 \sum_{k=0}^n & 2^k k^2 G_{p+k}H_{q+k} = \frac{1}{3} 2^{n+1} n^2 ( G_{p+n}H_{q+n+1} + G_{p+n-1}H_{q+n-2} ) \\
 & + \frac{1}{9} 2^{n+2} n ( 2G_{p+n}H_{q+n} - 6G_{p+n-1}H_{q+n-1} - 7G_{p+n-2}H_{q+n-2} ) \\
 & + \frac{2}{27} [ 2^n ( 70G_{p+n}H_{q+n} - 102G_{p+n-1}H_{q+n-1} - 173G_{p+n-2}H_{q+n-2} ) \\
 & \quad - 70G_pH_q + 102G_{p-1}H_{q-1} + 173G_{p-2}H_{q-2} ] \\
\end{split}
\end{equation}
\begin{equation}
\begin{split}
 \sum_{k=0}^n 2^k G_{p+2k}H_{q+2k} = \frac{1}{3} & [ 2^{n+1} ( 2G_{p+2n}H_{q+2n} - G_{p+2n-1}H_{q+2n-1} ) \\
 & + 2G_{p-1}H_{q-1} - G_pH_q ] \\
\end{split}
\end{equation}
\begin{equation}
\begin{split}
 \sum_{k=0}^n & 2^k G_{p+k}H_{q-k} \\
 = & \frac{1}{11} [ 2^{n+1} ( 6G_{p+n}H_{q-n} + 2G_{p+n-1}H_{q-n-1} + G_{p+n+1}H_{q-n+2} ) \\
 & \quad -G_pH_q - 4G_{p-1}H_{q-1} - 2G_{p+1}H_{q+2} ] \\
\end{split}
\end{equation}
Examples with $x=3$:
\begin{equation}
\begin{split}
 \sum_{k=0}^n 3^k G_{p+k}H_{q+k} = \frac{1}{4} & [ 3^{n+1} ( G_{p+n}H_{q+n+1} - G_{p+n-2}H_{q+n-3} ) \\
 & - G_{p-1}H_q + G_{p-3}H_{q-4} ] \\
\end{split}
\end{equation}
Examples with $x=1/2$:
\begin{equation}
\begin{split}
 \sum_{k=0}^n 2^{-k} G_{p+k}H_{q+k} = \frac{1}{3} & [ 2^{-n} ( 2G_{p+n+2}H_{q+n+2} - G_{p+n}H_{q+n} ) \\
 & - 2 ( 2G_{p+1}H_{q+1} - G_{p-1}H_{q-1} ) ] \\
\end{split}
\end{equation}
\begin{equation}
\begin{split}
 \sum_{k=0}^n & 2^{-k} k G_{p+k}H_{q+k} = \frac{1}{3} 2^{-n} n ( 2G_{p+n+2}H_{q+n+2} - G_{p+n}H_{q+n} ) \\
 & + \frac{2}{9} [ 2^{-n} ( 2G_{p+n}H_{q+n} - 6G_{p+n+1}H_{q+n+1} - 7G_{p+n+2}H_{q+n+2} ) \\
 & \quad - 2G_pH_q + 6G_{p+1}H_{q+1} + 7G_{p+2}H_{q+2} ] \\
\end{split}
\end{equation}
\begin{equation}
\begin{split}
 \sum_{k=0}^n & 2^{-k} k^2 G_{p+k}H_{q+k} = \frac{1}{3} 2^{-n} n^2 ( 2G_{p+n+2}H_{q+n+2} - G_{p+n}H_{q+n} ) \\
 & + \frac{4}{9} 2^{-n} n ( 2G_{p+n}H_{q+n} - 6G_{p+n+1}H_{q+n+1} - 7G_{p+n+2}H_{q+n+2} ) \\
 & - \frac{2}{27} [ 2^{-n} ( 70G_{p+n}H_{q+n} - 102G_{p+n+1}H_{q+n+1} - 173G_{p+n+2}H_{q+n+2} ) \\
 & \quad - 70G_pH_q + 102G_{p+1}H_{q+1} + 173G_{p+2}H_{q+2} ] \\
\end{split}
\end{equation}
\begin{equation}
\begin{split}
 \sum_{k=0}^n 2^{-k} G_{p+2k}H_{q+2k} = \frac{1}{3} & [ 2^{-n} ( 2G_{p+2n+1}H_{q+2n+1} - G_{p+2n}H_{q+2n} ) \\
 & + 4G_pH_q - 2G_{p+1}H_{q+1} ] \\
\end{split}
\end{equation}
\begin{equation}
\begin{split}
 \sum_{k=0}^n & 2^{-k} G_{p+k}H_{q-k} \\
 = & -\frac{1}{11} [ 2^{-n} ( G_{p+n}H_{q-n} + 4G_{p+n+1}H_{q-n+1} - 2G_{p+n-1}H_{q-n-2} ) \\
 & \quad - 12G_pH_q - 4G_{p+1}H_{q+1} + 2G_{p-1}H_{q-2} ] \\
\end{split}
\end{equation}
Examples with $x=1/3$:
\begin{equation}
\begin{split}
 \sum_{k=0}^n 3^{-k} G_{p+k}H_{q+k} = \frac{1}{4} & [ 3^{-n} ( G_{p+n+1}H_{q+n} - G_{p+n+3}H_{q+n+4} ) \\
 & - 3 ( G_pH_{q-1} - G_{p+2}H_{q+3} ) ] \\
\end{split}
\end{equation}
Examples with $x=-2$:
\begin{equation}
\begin{split}
 \sum_{k=0}^n & (-1)^k 2^k G_{p+k}H_{q+k} \\
 = \frac{1}{11} & [ (-1)^n 2^{n+1} ( 8G_{p+n}H_{q+n} - 2G_{p+n-1}H_{q+n-1} - G_{p+n+1}H_{q+n+1} ) \\
 & + 8G_{p-1}H_{q-1} - 2G_{p-2}H_{q-2} - G_pH_q ] \\
\end{split}
\end{equation}
\begin{equation}
\begin{split}
 \sum_{k=0}^n & (-1)^k 2^k G_{p+2k}H_{q+2k} \\
 = \frac{1}{57} & [ (-1)^n 2^{n+1} ( 18G_{p+2n}H_{q+2n} - 5G_{p+2n-1}H_{q+2n-1} + 4G_{p+2n+1}H_{q+2n+1} ) \\
 & + 21G_pH_q + 10G_{p-1}H_{q-1} - 8G_{p+1}H_{q+1} ] \\
\end{split}
\end{equation}
Examples with $x=-3$:
\begin{equation}
\begin{split}
 \sum_{k=0}^n & (-1)^k 3^k G_{p+k}H_{q+k} \\
 = \frac{1}{76} & [ (-1)^n 3^{n+1} ( 18G_{p+n}H_{q+n} - 14G_{p+n+1}H_{q+n+1} + 6G_{p+n+2}H_{q+n+2} ) \\
 & + 18G_{p-1}H_{q-1} - 14G_pH_q + 6G_{p+1}H_{q+1} ] \\
\end{split}
\end{equation}
Examples with $x=-1/2$:
\begin{equation}
\begin{split}
 \sum_{k=0}^n & (-1)^k 2^{-k} G_{p+k}H_{q+k} \\
 = \frac{1}{11} & [ (-1)^n 2^{-n} ( 4G_{p+n+1}H_{q+n+1} + 2G_{p+n-1}H_{q+n-1} - 5G_{p+n}H_{q+n} ) \\
 & + 2( 4G_pH_q + 2G_{p-2}H_{q-2} - 5G_{p-1}H_{q-1} ) ] \\
\end{split}
\end{equation}
\begin{equation}
\begin{split}
 \sum_{k=0}^n & (-1)^k 2^{-k} G_{p+2k}H_{q+2k} \\
 = \frac{1}{57} & [ (-1)^n 2^{-n} ( 5G_{p+2n}H_{q+2n} + 8G_{p+2n+2}H_{q+2n+2} - 6G_{p+2n+1}H_{q+2n+1} ) \\
 & + 2( 5G_{p-2}H_{q-2} + 8G_pH_q - 6G_{p-1}H_{q-1} ) ] \\
\end{split}
\end{equation}
Examples with $x=-1/3$:
\begin{equation}
\begin{split}
 \sum_{k=0}^n & (-1)^k 3^{-k} G_{p+k}H_{q+k} \\
 = \frac{1}{38} & [ (-1)^n 3^{-n} ( 9G_{p+n+1}H_{q+n+1} + 3G_{p+n-1}H_{q+n-1} - 7G_{p+n}H_{q+n} ) \\
 & + 3( 9G_pH_q + 3G_{p-2}H_{q-2} - 7G_{p-1}H_{q-1} ) ] \\
\end{split}
\end{equation}

\section{List of Cubic Examples}

\begin{equation}\label{cubic1}
\begin{split}
 \sum_{k=0}^n & G_{p+k}H_{q+k}K_{r+k} \\
 = \frac{1}{2} ( & G_{p+n-1}H_{q+n}K_{r+n} + G_{p+n}H_{q+n+1}K_{r+n+1} - G_{p+n+1}H_{q+n-1}K_{r+n-1} \\
  & - G_{p-2}H_{q-1}K_{r-1} - G_{p-1}H_qK_r + G_pH_{q-2}K_{r-2} ) \\
\end{split}
\end{equation}
\begin{equation}\label{cubic2}
\begin{split}
 \sum_{k=0}^n & (-1)^k G_{p+k}H_{q+k}K_{r+k} \\
 = \frac{1}{2} [ & (-1)^n ( G_{p+n+1}H_{q+n}K_{r+n} + G_{p+n}H_{q+n+2}K_{r+n+2} - G_{p+n+2}H_{q+n+1}K_{r+n+1} ) \\
  & + G_pH_{q-1}K_{r-1} + G_{p-1}H_{q+1}K_{r+1} - G_{p+1}H_qK_r ] \\
\end{split}
\end{equation}
\begin{equation}\label{cubic3}
\begin{split}
 \sum_{k=0}^n & G_{p+k}H_{q+k}K_{r-k} \\
 = \frac{1}{2} ( & G_{p+n+2}H_{q+n+2}K_{r-n} - G_{p+n}H_{q+n}K_{r-n-1} - G_{p+n+1}H_{q+n+1}K_{r-n-2} \\
  &  - G_{p+1}H_{q+1}K_{r+1} + G_{p-1}H_{q-1}K_r + G_pH_qK_{r-1} ) \\
\end{split}
\end{equation}
\begin{equation}
\begin{split}
 \sum_{k=0}^n & (-1)^k G_{p+k}H_{q+k}K_{r-k} \\
 = \frac{1}{2} [ & (-1)^n ( G_{p+n+1}H_{q+n+1}K_{r-n} + G_{p+n-1}H_{q+n-1}K_{r-n-1} - G_{p+n}H_{q+n}K_{r-n+1} ) \\
  &  + G_pH_qK_{r+1} + G_{p-2}H_{q-2}K_r - G_{p-1}H_{q-1}G_{r+2} ] \\
\end{split}
\end{equation}
\begin{equation}
\begin{split}
 \sum_{k=0}^n & G_{p+k}H_{q-k}K_{r-k} \\
 = \frac{1}{2} ( & G_{p+n+2}H_{q-n}K_{r-n} + G_{p+n-2}H_{q-n+1}K_{r-n+1} - G_{p+n}H_{q-n+2}K_{r-n+2} \\
  & - G_{p+1}H_{q+1}K_{r+1} - G_{p-3}H_{q+2}K_{r+2} + G_{p-1}H_{q+3}K_{r+3} ) \\
\end{split}
\end{equation}
\begin{equation}
\begin{split}
 \sum_{k=0}^n & (-1)^k G_{p+k}H_{q-k}K_{r-k} \\
 = \frac{1}{2} [ & (-1)^n ( G_{p+n+1}H_{q-n}K_{r-n} + G_{p+n}H_{q-n-2}K_{r-n-2} - G_{p+n+2}H_{q-n-1}K_{r-n-1} ) \\
  &  + G_pH_{q+1}K_{r+1} + G_{p-1}H_{q-1}K_{r-1} - G_{p+1}H_qK_r ] \\
\end{split}
\end{equation}
\begin{equation}
\begin{split}
 \sum_{k=0}^n & G_{p+2k}H_{q+2k}K_{r+2k} \\
 = \frac{1}{4} ( & G_{p+2n+2}H_{q+2n+2}K_{r+2n+2} - 3G_{p+2n+1}H_{q+2n+1}K_{r+2n+1} - G_{p+2n}H_{q+2n}K_{r+2n} \\
  & - G_pH_qK_r + 3G_{p-1}H_{q-1}K_{r-1} + G_{p-2}H_{q-2}K_{r-2} ) \\
\end{split}
\end{equation}
\begin{equation}
 \sum_{k=0}^n F_k^3 = \frac{1}{2} ( F_nF_{n+1}^2 - (-1)^n F_{n-1} + 1 )
\end{equation}
\begin{equation}
 \sum_{k=0}^n F_k^3 = \frac{1}{10} ( F_{3n+2} - 6(-1)^n F_{n-1} + 5 )
\end{equation}
\begin{equation}
 \sum_{k=0}^n k F_k^3 = \frac{1}{20} [ 2n F_{3n+2} - F_{3n+1} - 12(-1)^n ( n F_{n-1} + F_{n-3} ) + 25 ]
\end{equation}
\begin{equation}
\begin{split}
 \sum_{k=0}^n k^2 F_k^3 = \frac{1}{20} & [ 2n^2 F_{3n+2} - 2n F_{3n+1} + F_{3n+2} \\
  & - 12(-1)^n ( n^2 F_{n-1} + 2n F_{n-3} - F_{n-6} ) + 95 ] \\
\end{split}
\end{equation}
\begin{equation}
 \sum_{k=0}^n 2^{-k} F_k^3 = \frac{1}{25} \{ 2^{-n} [ F_{3n+2} + F_{3n+4} - 3(-1)^n ( F_{n-1} + F_{n+1} ) ] + 2 \}
\end{equation}
\begin{equation}
 \sum_{k=0}^n 3^{-k} F_k^3 = \frac{1}{110} \{ 3^{-n} [ 11 F_{3n+4} - 6(-1)^n ( 3 F_{n+1} - F_n ) ] - 15 \}
\end{equation}
\begin{equation}
 \sum_{k=0}^n F_{2k}^3 = \frac{1}{20} ( F_{6n+3} - 12 F_{2n+1} + 10 )
\end{equation}
\begin{equation}
 \sum_{k=0}^n (-1)^k F_k^3 = \frac{1}{2} ( (-1)^n F_n^2F_{n+1} - F_{n+2} + 1 )
\end{equation}
\begin{equation}
 \sum_{k=0}^n (-1)^k F_k^3 = \frac{1}{10} ( (-1)^n F_{3n+1} - 6 F_{n+2} + 5 )
\end{equation}
\begin{equation}
 \sum_{k=0}^n (-1)^k k F_k^3 = \frac{1}{20} [ (-1)^n ( 2n F_{3n+1} + F_{3n-1} ) - 12 ( n F_{n+2} - F_{n+3} ) - 25 ]
\end{equation}
\begin{equation}
\begin{split}
 \sum_{k=0}^n (-1)^k k^2 F_k^3 = \frac{1}{20} & [ (-1)^n ( 2n^2 F_{3n+1} + 2n F_{3n-1} - F_{3n-2} ) \\
  & - 12 ( n^2 F_{n+2} - 2n F_{n+3} + F_{n+6} ) + 95 ] \\
\end{split}
\end{equation}
\begin{equation}
 \sum_{k=0}^n L_k^3 = \frac{1}{2} ( L_nL_{n+1}^2 + 5(-1)^n L_{n-1} + 19 )
\end{equation}
\begin{equation}
 \sum_{k=0}^n L_k^3 = \frac{1}{2} ( L_{3n+2} + 6(-1)^n L_{n-1} + 19 )
\end{equation}
\begin{equation}
 \sum_{k=0}^n k L_k^3 = \frac{1}{4} [ 2n L_{3n+2} - L_{3n+1} + 12(-1)^n ( n L_{n-1} + L_{n-3} ) + 49 ]
\end{equation}
\begin{equation}
\begin{split}
 \sum_{k=0}^n k^2 L_k^3 = \frac{1}{4} & [ 2n^2 L_{3n+2} - 2n L_{3n+1} + L_{3n+2} \\
  & + 12(-1)^n ( n^2 L_{n-1} + 2n L_{n-3} - L_{n-6} ) + 213 ] \\
\end{split}
\end{equation}
\begin{equation}
 \sum_{k=0}^n F_k^4 = \frac{1}{125} [ L_{4n+4} - L_{4n} - 4(-1)^n ( L_{2n} + L_{2n+2} ) + 30 n + 15 ]
\end{equation}
\begin{equation}
 \sum_{k=0}^n L_k^4 = \frac{1}{5} [ L_{4n+4} - L_{4n} + 4(-1)^n ( L_{2n} + L_{2n+2} ) + 30 n + 55 ]
\end{equation}

\section{List of Quadratic Examples with Binomial\\ Coefficients}

Examples with $x=1$:
\begin{equation}\label{binomres1}
\begin{split}
 \sum_{k=0}^n \binom{n}{k} & G_{p+k}H_{q+k} = \frac{1}{5} \{ 5^{\lfloor n/2\rfloor} [ 2G_{p+n}(H_{q+1}-(-1)^nH_{q-1}) \\
 & \qquad + G_{p+n-1}(H_q-(-1)^nH_{q-2}) + G_{p+n+1}(H_{q+2}-(-1)^nH_q) ] \\
 & + \delta_{n,0} ( 3G_pH_q - G_{p-1}H_{q-1} - G_{p+1}H_{q+1} ) \} \\
\end{split}
\end{equation}
\begin{equation}\label{binomres2}
\begin{split}
 \sum_{k=0}^n & \binom{n}{k} k G_{p+k}H_{q+k} = \frac{1}{5} \{ 5^{\lfloor (n-1)/2\rfloor} n [ 2G_{p+n}(H_{q+2}+(-1)^nH_q) \\
 & \qquad + G_{p+n-1}(H_{q+1}+(-1)^nH_{q-1}) + G_{p+n+1}(H_{q+3}+(-1)^nH_{q+1}) ] \\
 & - \delta_{n,1} ( 3G_pH_q - G_{p-1}H_{q-1} - G_{p+1}H_{q+1} ) \} \\
\end{split}
\end{equation}
\begin{equation}\label{binomres3}
\begin{split}
 \sum_{k=0}^n & \binom{n}{k} k^2 G_{p+k}H_{q+k} = \frac{1}{5} \{ 5^{\lfloor (n-1)/2\rfloor} n [ 2G_{p+n}(H_{q+2}+(-1)^nH_q) \\
 & \qquad + G_{p+n-1}(H_{q+1}+(-1)^nH_{q-1}) + G_{p+n+1}(H_{q+3}+(-1)^nH_{q+1}) ] \\
 & + 5^{\lfloor n/2\rfloor-1} n(n-1) [ 2G_{p+n}(H_{q+3}-(-1)^nH_{q+1}) \\
 & \qquad + G_{p+n-1}(H_{q+2}-(-1)^nH_q) + G_{p+n+1}(H_{q+4}-(-1)^nH_{q+2}) ] \\
 & + ( 2\delta_{n,2} - \delta_{n,1} ) ( 3G_pH_q - G_{p-1}H_{q-1} - G_{p+1}H_{q+1} ) \} \\
\end{split}
\end{equation}
\begin{equation}\label{binomres3}
\begin{split}
 \sum_{k=0}^n & \binom{n}{k} k^3 G_{p+k}H_{q+k} = \frac{1}{5} \{ 5^{\lfloor (n-1)/2\rfloor-1} n [ 5 ( 2G_{p+n}(H_{q+2}+(-1)^nH_q) \\
 & \qquad + G_{p+n-1}(H_{q+1}+(-1)^nH_{q-1}) + G_{p+n+1}(H_{q+3}+(-1)^nH_{q+1}) ) \\
 & \quad + (n-1)(n-2) (  2G_{p+n}(H_{q+4}+(-1)^nH_{q+2}) \\
 & \qquad + G_{p+n-1}(H_{q+3}+(-1)^nH_{q+1}) + G_{p+n+1}(H_{q+5}+(-1)^nH_{q+3}) ) ] \\
 & + 5^{\lfloor n/2\rfloor-1} 3n(n-1) [ 2G_{p+n}(H_{q+3}-(-1)^nH_{q+1}) \\
 & \qquad + G_{p+n-1}(H_{q+2}-(-1)^nH_q) + G_{p+n+1}(H_{q+4}-(-1)^nH_{q+2}) ] \\
 & + ( 6\delta_{n,2} - \delta_{n,1} - 6\delta_{n,3} ) ( 3G_pH_q - G_{p-1}H_{q-1} - G_{p+1}H_{q+1} ) \} \\
\end{split}
\end{equation}
\begin{equation}
\begin{split}
 \sum_{k=0}^n \binom{n}{k} & G_{p+2k}H_{q+k} \\
  = \frac{1}{5} & [ 2^n ( 2G_{p+n}H_{q+n} + G_{p+n-1}H_{q+n-1} + G_{p+n+1}H_{q+n+1} ) \\
  & + 3G_pH_{q+n} - G_{p-1}H_{q+n-1} - G_{p+1}H_{q+n+1} ] \\
\end{split}
\end{equation}
\begin{equation}
\begin{split}
 \sum_{k=0}^n \binom{n}{k} & k G_{p+2k}H_{q+k} \\
  = \frac{1}{5} & n [ 2^{n-1} ( 2G_{p+n+1}H_{q+n} + G_{p+n}H_{q+n-1} + G_{p+n+2}H_{q+n+1} ) \\
  & + 3G_{p+2}H_{q+n} - G_{p+1}H_{q+n-1} - G_{p+3}H_{q+n+1} ] \\
\end{split}
\end{equation}\begin{equation}
\begin{split}
 \sum_{k=0}^n \binom{n}{k} & k^2 G_{p+2k}H_{q+k} \\
  = \frac{1}{5} & n \{ 2^{n-2} [ n ( 2G_{p+n+2}H_{q+n} + G_{p+n+1}H_{q+n-1} + G_{p+n+3}H_{q+n+1} ) \\
 & \qquad + 2G_{p+n-1}H_{q+n} + G_{p+n-2}H_{q+n-1} + G_{p+n}H_{q+n+1} ] \\ 
  & + n ( 3G_{p+4}H_{q+n} - G_{p+3}H_{q+n-1} - G_{p+5}H_{q+n+1} )  \\
 & \qquad - 3G_{p+3}H_{q+n} + G_{p+2}H_{q+n-1} + G_{p+4}H_{q+n+1} \} \\
\end{split}
\end{equation}
\begin{equation}
\begin{split}
 \sum_{k=0}^n \binom{n}{k} & G_{p+2k}H_{q+2k} \\
  = \frac{1}{5} & [ 3^n ( 2G_{p+n}H_{q+n} + G_{p+n-1}H_{q+n-1} + G_{p+n+1}H_{q+n+1} ) \\
  & +  2^n ( 3G_pH_q - G_{p-1}H_{q-1} - G_{p+1}H_{q+1} )  ] \\
\end{split}
\end{equation}
\begin{equation}
\begin{split}
 \sum_{k=0}^n \binom{n}{k} & k G_{p+2k}H_{q+2k} \\
  = \frac{1}{5} & n [ 3^{n-1} ( 2G_{p+n+1}H_{q+n+1} + G_{p+n}H_{q+n} + G_{p+n+2}H_{q+n+2} ) \\
  & +  2^{n-1} ( 3G_pH_q - G_{p-1}H_{q-1} - G_{p+1}H_{q+1} )  ] \\
\end{split}
\end{equation}
\begin{equation}
\begin{split}
 \sum_{k=0}^n \binom{n}{k} & k^2 G_{p+2k}H_{q+2k} \\
  = \frac{1}{5} & n \{ 3^{n-2} [ n ( 2G_{p+n+2}H_{q+n+2} + G_{p+n+1}H_{q+n+1} + G_{p+n+3}H_{q+n+3} ) \\
 & \qquad + 2G_{p+n}H_{q+n} + G_{p+n-1}H_{q+n-1} + G_{p+n+1}H_{q+n+1} ] \\
  & +  2^{n-2} (n+1) ( 3G_pH_q - G_{p-1}H_{q-1} - G_{p+1}H_{q+1} )  \} \\
\end{split}
\end{equation}
\begin{equation}
\begin{split}
 \sum_{k=0}^n \binom{n}{k} G_{p+k}H_{q-k} 
  = & \frac{1}{5} [ 2^n ( 2G_pH_q + G_{p-1}H_{q-1} + G_{p+1}H_{q+1} ) \\
  & + 3G_pH_{q-n} - G_{p-1}H_{q-n-1} - G_{p+1}H_{q-n+1} ] \\
\end{split}
\end{equation}
\begin{equation}
\begin{split}
 \sum_{k=0}^n \binom{n}{k} k G_{p+k}H_{q-k} 
  = & \frac{1}{5} n [ 2^{n-1} ( 2G_pH_q + G_{p-1}H_{q-1} + G_{p+1}H_{q+1} ) \\
  & + 3G_{p+1}H_{q-n} - G_pH_{q-n-1} - G_{p+2}H_{q-n+1} ] \\
\end{split}
\end{equation}
\begin{equation}
\begin{split}
 \sum_{k=0}^n \binom{n}{k} k^2 G_{p+k}H_{q-k} 
  = & \frac{1}{5} n [ 2^{n-2} (n+1) ( 2G_pH_q + G_{p-1}H_{q-1} + G_{p+1}H_{q+1} ) \\
  & + n ( 3G_{p+2}H_{q-n} - G_{p+1}H_{q-n-1} - G_{p+3}H_{q-n+1} ) \\
 &  - 3G_pH_{q-n} + G_{p-1}H_{q-n-1} + G_{p+1}H_{q-n+1} ] \\
\end{split}
\end{equation}
\begin{equation}
\begin{split}
 \sum_{k=0}^n & \binom{n}{k} k^3 G_{p+k}H_{q-k} \\
  = & \frac{1}{5} n [ 2^{n-3} n(n+3) ( 2G_pH_q + G_{p-1}H_{q-1} + G_{p+1}H_{q+1} ) \\
  & + n^2 ( 3G_{p+3}H_{q-n} - G_{p+2}H_{q-n-1} - G_{p+4}H_{q-n+1} ) \\
  & -(3n-2) ( 3G_{p+1}H_{q-n} - G_pH_{q-n-1} - G_{p+2}H_{q-n+1} ) \\
 &  - 3G_pH_{q-n} + G_{p-1}H_{q-n-1} + G_{p+1}H_{q-n+1} ] \\
\end{split}
\end{equation}
\begin{equation}
\begin{split}
 \sum_{k=0}^n \binom{n}{k} G_{p+2k}H_{q-k} 
  = & \frac{1}{5} [ 2G_{p+n}H_{q+n} + G_{p+n-1}H_{q+n-1} + G_{p+n+1}H_{q+n+1} \\
  & + 2^n ( 3G_pH_{q-n} - G_{p-1}H_{q-n-1} - G_{p+1}H_{q-n+1} ) ] \\
\end{split}
\end{equation}
\begin{equation}
\begin{split}
 \sum_{k=0}^n \binom{n}{k} G_{p+2k}H_{q-2k} 
  = & \frac{1}{5} [ 2^n( 2G_pH_q + G_{p-1}H_{q-1} + G_{p+1}H_{q+1} ) \\
  & + 3^n ( 3G_{p+2n}H_q - G_{p+2n-1}H_{q-1} - G_{p+2n+1}H_{q+1} ) ] \\
\end{split}
\end{equation}
\begin{equation}
\begin{split}
 \sum_{k=0}^n \binom{n}{k} k & G_{p+2k}H_{q-2k} 
  = \frac{1}{5} n [ 2^{n-1}( 2G_pH_q + G_{p-1}H_{q-1} + G_{p+1}H_{q+1} ) \\
  & + 3^{n-1} ( 3G_{p+2n+2}H_q - G_{p+2n+1}H_{q-1} - G_{p+2n+3}H_{q+1} ) ] \\
\end{split}
\end{equation}
\begin{equation}
\begin{split}
 \sum_{k=0}^n & \binom{n}{k} k^2 G_{p+2k}H_{q-2k} \\
  & = \frac{1}{5} n \{ 2^{n-2}(n+1)( 2G_pH_q + G_{p-1}H_{q-1} + G_{p+1}H_{q+1} ) \\
  & \qquad + 3^{n-2} [ n ( 3G_{p+2n+4}H_q - G_{p+2n+3}H_{q-1} - G_{p+2n+5}H_{q+1} )  \\
 & \qquad + 3G_{p+2n}H_q - G_{p+2n-1}H_{q-1} - G_{p+2n+1}H_{q+1} ] \} \\
\end{split}
\end{equation}
Examples with $x=-1$:
\begin{equation}
\begin{split}
 \sum_{k=0}^n \binom{n}{k} & (-1)^k G_{p+k}H_{q+k} \\
  = \frac{1}{5} & [ (-1)^n ( 2G_{p+n}H_q + G_{p+n-1}H_{q-1} + G_{p+n+1}H_{q+1} ) \\
  & +  2^n ( 3G_pH_q - G_{p-1}H_{q-1} - G_{p+1}H_{q+1} )  ] \\
\end{split}
\end{equation}
\begin{equation}
\begin{split}
 \sum_{k=0}^n \binom{n}{k} & (-1)^k k G_{p+k}H_{q+k} \\
  = \frac{1}{5} & n [ (-1)^n ( 2G_{p+n+1}H_q + G_{p+n}H_{q-1} + G_{p+n+2}H_{q+1} ) \\
  & +  2^{n-1} ( 3G_pH_q - G_{p-1}H_{q-1} - G_{p+1}H_{q+1} )  ] \\
\end{split}
\end{equation}
\begin{equation}
\begin{split}
 \sum_{k=0}^n \binom{n}{k} & (-1)^k k^2 G_{p+k}H_{q+k} \\
  = \frac{1}{5} & n \{ (-1)^n [ n ( 2G_{p+n+2}H_q + G_{p+n+1}H_{q-1} + G_{p+n+3}H_{q+1} ) \\
  & \qquad - 2G_{p+n}H_q - G_{p+n-1}H_{q-1} - G_{p+n+1}H_{q+1} ] \\
  & +  2^{n-2}(n+1) ( 3G_pH_q - G_{p-1}H_{q-1} - G_{p+1}H_{q+1} )  \} \\
\end{split}
\end{equation}
\begin{equation}
\begin{split}
 \sum_{k=0}^n \binom{n}{k} & (-1)^k k^3 G_{p+k}H_{q+k} \\
  = \frac{1}{5} & n \{ (-1)^n [ n^2 ( 2G_{p+n+3}H_q + G_{p+n+2}H_{q-1} + G_{p+n+4}H_{q+1} ) \\
  & \qquad -(3n-2) ( 2G_{p+n+1}H_q + G_{p+n}H_{q-1} + G_{p+n+2}H_{q+1} ) \\
  & \qquad - 2G_{p+n}H_q - G_{p+n-1}H_{q-1} - G_{p+n+1}H_{q+1} ] \\
  & +  2^{n-3}n(n+3) ( 3G_pH_q - G_{p-1}H_{q-1} - G_{p+1}H_{q+1} )  \} \\
\end{split}
\end{equation}
\begin{equation}
\begin{split}
 \sum_{k=0}^n \binom{n}{k} & (-1)^k G_{p+2k}H_{q+k} \\
  = \frac{1}{5} & [ (-1)^n 2^n( 2G_{p+n}H_q + G_{p+n-1}H_{q-1} + G_{p+n+1}H_{q+1} ) \\
  & +  3G_{p+2n}H_q - G_{p+2n-1}H_{q-1} - G_{p+2n+1}H_{q+1}  ] \\
\end{split}
\end{equation}
\begin{equation}
\begin{split}
 \sum_{k=0}^n \binom{n}{k} & (-1)^k k G_{p+2k}H_{q+k} \\
  = \frac{1}{5} & n [ (-1)^n 2^{n-1}( 2G_{p+n+2}H_q + G_{p+n+1}H_{q-1} + G_{p+n+3}H_{q+1} ) \\
  & +  3G_{p+2n-1}H_q - G_{p+2n-2}H_{q-1} - G_{p+2n}H_{q+1}  ] \\
\end{split}
\end{equation}
\begin{equation}
\begin{split}
 \sum_{k=0}^n \binom{n}{k} & (-1)^k k^2 G_{p+2k}H_{q+k} \\
  = \frac{1}{5} & n \{ (-1)^n 2^{n-2}[ n(2G_{p+n+4}H_q + G_{p+n+3}H_{q-1} + G_{p+n+5}H_{q+1} ) \\
 & \qquad - 2G_{p+n+1}H_q - G_{p+n}H_{q-1} - G_{p+n+2}H_{q+1} ] \\
  & + n ( 3G_{p+2n-2}H_q - G_{p+2n-3}H_{q-1} - G_{p+2n-1}H_{q+1}  ) \\
 & +  3G_{p+2n-3}H_q - G_{p+2n-4}H_{q-1} - G_{p+2n-2}H_{q+1} \} \\
\end{split}
\end{equation}
\begin{equation}
\begin{split}
 & \sum_{k=0}^n \binom{n}{k} (-1)^k G_{p+2k}H_{q+2k} \\
  & = \frac{1}{5} \{ (-1)^n 5^{\lfloor n/2\rfloor} [ 2G_{p+n}(H_{q+n+1}-(-1)^nH_{q+n-1}) + \\ 
 & \qquad+ G_{p+n-1}(H_{q+n}-(-1)^nH_{q+n-2}) + G_{p+n+1}(H_{q+n+2}-(-1)^nH_{q+n}) ] \\
  & \quad + \delta_{n,0} ( 3G_pH_q - G_{p-1}H_{q-1} - G_{p+1}H_{q+1} )  \} \\
\end{split}
\end{equation}
\begin{equation}
\begin{split}
 \sum_{k=0}^n & \binom{n}{k} (-1)^k G_{p+k}H_{q-k} \\
 & = \frac{1}{5} \{  5^{\lfloor n/2\rfloor} [ 3(G_{p+n+1}-(-1)^nG_{p+n-1})H_q \\
 & \qquad - (G_{p+n}-(-1)^nG_{p+n-2}) H_{q-1} - (G_{p+n+2}-(-1)^nG_{p+n})H_{q+1}  ] \\
 & \quad + \delta_{n,0} ( 2G_pH_q + G_{p-1}H_{q-1} + G_{p+1}H_{q+1} ) \} \\
\end{split}
\end{equation}
\begin{equation}
\begin{split}
 & \sum_{k=0}^n \binom{n}{k} (-1)^k k G_{p+k}H_{q-k} \\
 & = \frac{1}{5} \{  5^{\lfloor (n-1)/2\rfloor} n [ 3(G_{p+n+2}+(-1)^nG_{p+n})H_q \\
 & \qquad - (G_{p+n+1}+(-1)^nG_{p+n-1}) H_{q-1} - (G_{p+n+3}+(-1)^nG_{p+n+1})H_{q+1}  ] \\
 & \quad - \delta_{n,1} ( 2G_pH_q + G_{p-1}H_{q-1} + G_{p+1}H_{q+1} ) \} \\
\end{split}
\end{equation}
\begin{equation}
\begin{split}
 & \sum_{k=0}^n \binom{n}{k} (-1)^k k^2 G_{p+k}H_{q-k} \\
 & = \frac{1}{5} \{  5^{\lfloor (n-1)/2\rfloor} n [ 3(G_{p+n+2}+(-1)^nG_{p+n})H_q \\
 & \qquad - (G_{p+n+1}+(-1)^nG_{p+n-1}) H_{q-1} - (G_{p+n+3}+(-1)^nG_{p+n+1})H_{q+1}  ] \\
 & \quad +  5^{\lfloor n/2\rfloor-1} n(n-1) [ 3(G_{p+n+3}-(-1)^nG_{p+n+1})H_q \\
 & \qquad - (G_{p+n+2}-(-1)^nG_{p+n}) H_{q-1} - (G_{p+n+4}-(-1)^nG_{p+n+2})H_{q+1}  ] \\
 & \quad +( 2\delta_{n,2} - \delta_{n,1} ) ( 2G_pH_q + G_{p-1}H_{q-1} + G_{p+1}H_{q+1} ) \} \\
\end{split}
\end{equation}
\begin{equation}
\begin{split}
 \sum_{k=0}^n & \binom{n}{k} (-1)^k G_{p+2k}H_{q-k} \\
 & = \frac{1}{5} [ (-1)^n ( 2G_{p-n}H_q + G_{p-n-1}H_{q-1} + G_{p-n+1}H_{q+1} ) \\
 & \quad + 2^n ( 3G_{p+2n}H_q - G_{p+2n-1}H_{q-1} - G_{p+2n+1}H_{q+1} ) ] \\
\end{split}
\end{equation}
\begin{equation}
\begin{split}
 & \sum_{k=0}^n \binom{n}{k} (-1)^k G_{p+2k}H_{q-2k} \\
 & = \frac{1}{5} \{  5^{\lfloor n/2\rfloor} [ 3(G_{p+n+1}-(-1)^nG_{p+n-1})H_{q-n} \\
 & \qquad - (G_{p+n}-(-1)^nG_{p+n-2}) H_{q-n-1} - (G_{p+n+2}-(-1)^nG_{p+n})H_{q-n+1}  ] \\
 & \quad + \delta_{n,0} ( 2G_pH_q + G_{p-1}H_{q-1} + G_{p+1}H_{q+1} ) \} \\
\end{split}
\end{equation}

\section{List of Cubic Examples with Binomial\\ Coefficients}

\begin{equation}
\begin{split}
 & \sum_{k=0}^n \binom{n}{k} G_{p+k}H_{q+k}K_{r+k} \\
 & = \frac{1}{5} [ 2^n ( G_pH_{q+n}K_{r+n} + G_{p+1}H_{q+n+1}K_{r+n+1} - G_{p-1}H_{q+n-1}K_{r+n-1} ) \\
 & \quad + G_pH_qK_{r+n} + G_{p-1}H_{q-1}K_{r+n+1} - G_{p+1}H_{q+1}K_{r+n-1} \\
 & \quad + (-1)^n ( 3G_pH_q - G_{p-1}H_{q-1} - G_{p+1}H_{q+1} ) K_{r-n} ] \\
\end{split}
\end{equation}
\begin{equation}
\begin{split}
 & \sum_{k=0}^n \binom{n}{k} (-1)^k G_{p+k}H_{q+k}K_{r+k} \\
 & = \frac{1}{5} [ (-1)^n 2^n ( G_pH_{q+n}K_r + G_{p+1}H_{q+n+1}K_{r+1} - G_{p-1}H_{q+n-1}K_{r-1} ) \\
 & \quad + G_pH_{q+2n}K_r + G_{p-1}H_{q+2n-1}K_{r+1} - G_{p+1}H_{q+2n+1}K_{r-1} \\
 & \quad + ( 3G_pH_q - G_{p-1}H_{q-1} - G_{p+1}H_{q+1} ) K_{r+2n} ] \\
\end{split}
\end{equation}
\begin{equation}
\begin{split}
 & \sum_{k=0}^n \binom{n}{k} G_{p+k}H_{q+k}K_{r-k} \\
 & = \frac{1}{5} [ 2^n ( G_pH_qK_{r-n} + G_{p-1}H_{q-1}K_{r-n+1} - G_{p+1}H_{q+1}K_{r-n-1} ) \\
 & \quad + G_pH_{q+n}K_{r+n} + G_{p+1}H_{q+n+1}K_{r+n+1} - G_{p-1}H_{q+n-1}K_{r+n-1} \\
 & \quad + ( 3G_pH_q - G_{p-1}H_{q-1} - G_{p+1}H_{q+1} ) K_{r-2n} ] \\
\end{split}
\end{equation}
\begin{equation}
\begin{split}
 & \sum_{k=0}^n \binom{n}{k} (-1)^k G_{p+k}H_{q+k}K_{r-k} \\
 & = \frac{1}{5} [ 2^n ( G_{p+n}H_{q+n}K_r + G_{p+n-1}H_{q+n-1}K_{r+1} - G_{p+n+1}H_{q+n+1}K_{r-1} ) \\
 & \quad + (-1)^n ( G_pH_{q-n}K_r + G_{p+1}H_{q-n+1}K_{r+1} - G_{p-1}H_{q-n-1}K_{r-1} ) \\
 & \quad + ( 3G_pH_q - G_{p-1}H_{q-1} - G_{p+1}H_{q+1} ) K_{r+n} ] \\
\end{split}
\end{equation}
\begin{equation}
\begin{split}
 & \sum_{k=0}^n \binom{n}{k} G_{p+k}H_{q-k}K_{r-k} \\
 & = \frac{1}{5} [ 2^n ( G_{p+2n}H_qK_r + G_{p+2n+1}H_{q-1}K_{r-1} - G_{p+2n-1}H_{q+1}K_{r+1} ) \\
 & \quad + G_pH_{q+n}K_r + G_{p-1}H_{q+n+1}K_{r-1} - G_{p+1}H_{q+n-1}K_{r+1} \\
 & \quad + ( 2G_pH_q + G_{p-1}H_{q-1} + G_{p+1}H_{q+1} ) K_{r+n} ] \\
\end{split}
\end{equation}
\begin{equation}
\begin{split}
 & \sum_{k=0}^n \binom{n}{k} (-1)^k G_{p+k}H_{q-k}K_{r-k} \\
 & = \frac{1}{5} [ 2^n ( G_pH_qK_{r-n} + G_{p+1}H_{q-1}K_{r-n-1} - G_{p-1}H_{q+1}K_{r-n+1} ) \\
 & \quad + G_{p+n}H_qK_{r+n} + G_{p+n-1}H_{q+1}K_{r+n-1} - G_{p+n+1}H_{q-1}K_{r+n+1} \\
 & \quad + ( 2G_pH_q + G_{p-1}H_{q-1} + G_{p+1}H_{q+1} ) K_{r-2n} ] \\
\end{split}
\end{equation}

\pdfbookmark[0]{References}{}

\end{document}